\documentclass[12pt]{article}

\usepackage{amsmath,amsthm,amsfonts,color,hyperref,anysize,amssymb,enumitem,sectsty,color,dsfont}
\usepackage[section]{placeins}

\marginsize{0.8in}{0.8in}{0.8in}{0.8in}

\newcommand{\ep}{\epsilon}
\newcommand{\rr}{\mathbb{R}}

\newcommand{\ra}{\rightarrow}

\newcommand{\pc}{p^0}
\newcommand{\ps}{p^*}
\newcommand{\tp}{\tilde{p}}
\newcommand{\tP}{\tilde{P}}
\newcommand{\comm}[1]{\qquad\mbox{(#1)}}
\newcommand{\Span}{\mbox{span}}
\newcommand{\eF}{\epsilon_{\text{\tiny F}}}
\newcommand{\eB}{\epsilon_{\text{\tiny B}}}
\newcommand{\PP}[4]{\tP_{#1}^{[#2,#3]}\left(#4\right)}

\newcommand{\hide}[1]{}
\hide{
     \newcommand{\qed}{\nobreak \ifvmode \relax \else
     \ifdim\lastskip<1.5em \hskip-\lastskip
     \hskip1.5em plus0em minus0.5em \fi \nobreak
     \vrule height0.75em width0.5em depth0.25em\fi}
}

\newtheorem{theorem}{Theorem}
\newtheorem{lemma}[theorem]{Lemma}

\newtheorem{cor}[theorem]{Corollary}

\newtheorem{defn}{Definition}

\newenvironment{pf}{\begin{proof}[\emph{\textbf{Proof: }}]}{\end{proof}}
\newenvironment{pfof}[1]{\begin{proof}[\emph{\textbf{Proof of #1: }}]}{\end{proof}}

\title{Amortized Analysis on Asynchronous Gradient Descent}

\author{Yun Kuen Cheung\thanks{Most of the work done while at Courant Institute, NYU.}\\University of Vienna
\and
Richard Cole\\Courant Institute, NYU
}
\date{}

\begin{document}\setlength{\parindent}{0.2in}\setlength{\parskip}{0.04in}
\maketitle

\begin{abstract}
\thispagestyle{empty}
\setcounter{page}{0}
Gradient descent is an important class of iterative algorithms for minimizing convex functions.
Classically, gradient descent has been a sequential and synchronous process.
Distributed and asynchronous variants of gradient descent have been studied since the 1980s,
and they have been experiencing a resurgence due to demand from large-scale machine learning problems running on multi-core processors.

We provide a version of asynchronous gradient descent (AGD) in which communication between cores is minimal
and for which there is little synchronization overhead.
We also propose a new timing model for its analysis.
With this model, we give the first amortized analysis of AGD on convex functions.
The amortization allows for bad updates (updates that increase the value of the convex function);
in contrast, most prior work makes the strong assumption that every update must be significantly improving.

Typically, the step sizes used in AGD are smaller than those used in its synchronous counterpart.
We provide a method to determine the step sizes in AGD based on the Hessian entries for the convex function.
In certain circumstances, the resulting step sizes are a constant fraction of those used in the corresponding synchronous algorithm,
enabling the overall performance of AGD to improve linearly with the number of cores.

We give two applications of our amortized analysis:
\begin{itemize}
\item We show that our AGD algorithm can be applied to two classes of problems
which have huge problem sizes in applications and consequently can benefit substantially from parallelism.
The first class of problems is to solve linear systems $Ap=b$,
where the $A$ are symmetric and positive definite matrices.
The second class of problems is to minimize convex functions of the form
$\sum_{i=1}^n f_i(p_i) + \frac{1}{2} \|Ap-b\|^2$, where the $f_i$ are convex differentiable univariate functions.

\item We show that a version of asynchronous tatonnement, a simple distributed price update dynamic,
converges toward the market equilibrium in Fisher markets with buyers having complementary-CES or Leontief utility functions.
\end{itemize}
\end{abstract} 

\newpage
\thispagestyle{plain}
\setcounter{page}{1}

\section{Introduction}

Gradient descent, an important class of iterative algorithms for minimizing convex functions,
is a key subroutine in many computational problems.
Broadly speaking, gradient descent proceeds by iteratively moving
in the direction of the negative gradient of the convex function.
Classically, gradient descent is a sequential and synchronous process.
Distributed and asynchronous variants have also been studied,
starting with the work of Tsitsiklis et al.~\cite{TBA1986} in the 1980s;
more recent results include~\cite{BT2000,Borkar1998}.
Distributed and asynchronous gradient descent has been experiencing a resurgence of attention,
particularly in computational learning theory~\cite{LSZ2009,NRRW2011},
due to recent advances in multi-core parallel processing technology
and a strong demand for speeding-up large-scale gradient descent problems via parallelism.

Gradient descent proceeds by repeatedly updating the coordinates of the argument to the convex function.
A few key common issues arise in any distributed and asynchronous iterative implementation
and their improper handling may lead to performance-destroying overhead costs.
\begin{itemize}
\item In some implementations (e.g.~\cite{NRRW2011}), different
cores\footnote{These observations apply to any multi-processor system.}
may update the same component.
Without proper coordination, the progress made by one core can be overwritten,
and if such overwriting persists, in the worst case the system can fail to reach the desired result.

This difficulty can be avoided by block component descent -- each coordinate is updated by exactly one core.
This is the approach we use in our Asynchronous Gradient Descent (AGD) algorithm.
The approach has been used previously in a round-robin manner~\cite{LSZ2009},
but our AGD algorithm does not require the updates to proceed in any particular order.

\item The cores need to follow a communication protocol in order to communicate/broadcast their updates.
Communication is often relatively slow compared to computation,
so reducing the need for communication can lead to a significant improvement in system performance.
Also, when there is delay in communication, cores may use outdated information for the next update,
which is a critical issue for asynchronous systems.

One common approach is to assume that the system has \emph{bounded asynchrony},
i.e.~the delay in communication is bounded by a positive constant.
Typically, there is a need to wait for updates from the other cores, and the bounded asynchrony simply bounds the waiting time.
We will use the bounded asynchrony assumption, but our AGD algorithm will have \emph{no waiting}:
updates will always be based on the information at hand;
bounded asynchrony just guarantees that it is not too dated.

\item Often, the computation of one core needs the results computed by another core,
implying the computations of the different cores must be in a correct order to ensure correctness and to reduce core waiting time.
Typically this is achieved via a synchronization protocol, which often requires that all cores follow a global clock.
However, such protocols can be costly and even impractical in some circumstances.

As we shall see, our AGD algorithm needs essentially no synchronization apart from
an initial synchronization to align the starting times of all cores.
\end{itemize}

Broadly speaking, most prior work follows the asynchrony model proposed in~\cite{TBA1986}, in which time is discretized.
Our AGD algorithm allows each core to proceed at its own pace.
This allows for varying loads, for different updates having varied costs, for interruptions,
and more generally for variations in the completion times of updates.
To support this, in our model, time is continuous.
To ensure progress, we require that each component be updated at least once in each time unit,
but do not impose an upper bound on the frequency of updates.
A more formal description of our model will be given in Section~\ref{sect:model}.

We consider a robust family of AGD algorithms, and using our timing model,
we give a new amortized analysis which shows each algorithm converges to the minimal value of the underlying function.
Most prior work made the strong assumption that each update yields a significant improvement.
Our analysis, however, allows for bad individual updates (updates that increase the value of the convex function),
which seem to be unavoidable in general.
In our AGD algorithm, every update leads to errors in subsequent gradient measurements at other cores.
A natural question to ask is whether such errors can propagate and be persistent and whether they might, in the worst case,
prohibit convergence toward a minimal point.
Our amortized analysis shows that this will not happen when the step sizes used in the AGD algorithm are suitably bounded.
The following observation forms a key part of the analysis:
if there is a bad update to one component,
it can only be due to some recent good updates to other components, or to chaining of this effect.
We use a carefully designed potential function, which \emph{saves}
a portion of the gains due to good updates, to pay for the bad updates.
The amortized analysis will be presented in Section~\ref{sect:analysis}.

Typically the step sizes used in AGD are smaller than those used in its synchronous counterpart.
Our AGD algorithm determines the step sizes based on the Hessian of the underlying function.
In certain circumstances, the step sizes in our AGD can be a constant fraction of those used in its synchronous counterpart,
ensuring that the number of rounds of updates performed by the AGD algorithm is
within a constant of the analogous upper bound for the synchronous version.
Note that AGD avoids the synchronization costs of its synchronous counterpart, which are a practical concern~\cite{NRRW2011}.

\paragraph{Application: Solving Matrix Systems in Parallel}
We begin by considering two problems in which bad updates are possible in an asynchronous setting.
A linear system is the problem of finding $p\in\rr^n$ that satisfies $Ap = b$, where $A\in\rr^{m\times n}$ and $b\in \rr^m$ are the inputs.
As is well-known, if $A$ is a symmetric and positive definite matrix,
solving the linear system is equivalent to finding the minimum point of a strongly convex function,
so our AGD algorithm can be applied.

Nesterov \cite{Nesterov2012} discusses the following class of optimization problems:
minimizing convex functions of the form $\sum_{i=1}^n f_i(p_i) + \frac{1}{2} \|Ap-b\|^2$,
where the $f_i$ are convex differentiable univariate functions.
The size of such problems can be huge in practice, and input/data can be distributed in space and time,
so time synchronization is costly and even impractical.
One important feature of our AGD algorithm is to allow the use of data that are \emph{variously dated}.
As we will see, this hugely reduces the need for synchronization.
More details are given in Section~\ref{sect:linear-systems}.

\paragraph{Application: Asynchronous Tatonnement in Fisher Markets}
We show that an asynchronous tatonnement converges toward the market equilibrium in two classes of Fisher markets.

The concept of a market equilibrium was first proposed by Walras \cite{Walras1874}.
Walras also proposed an algorithmic approach for finding equilibrium prices, namely to
adjust prices by tatonnement: upward if there is too much demand and downward if too little.
Since then, the study of market equilibria and tatonnement have received much attention in economics, operations research,
and most recently in computer science~\cite{ABH1959,Uzawa1960,CMV2005,PY2010}.
Underlying many of these works is the issue of what are plausible price adjustment mechanisms
and in what types of markets they attain a market equilibrium.

The tatonnements studied in prior work have mostly been continuous, or discrete and synchronous.
Observing that real-world market dynamics are highly distributed and hence presumably asynchronous,
Cole and Fleischer~\cite{CF2008} initiated the study of asynchronous tatonnement with their \emph{Ongoing market model},
a market model incorporating update dynamics.

Cheung, Cole and Devanur~\cite{CCD2013} showed that tatonnement is equivalent to gradient descent
on a convex function for several classes of Fisher markets, and consequently
that a suitable synchronous tatonnement converges toward the market equilibrium
in two classes of markets: complementary-CES Fisher markets and Leontief Fisher markets.
This equivalence also enables us to apply our amortized analysis
to show that the corresponding asynchronous version of tatonnement
converges toward the market equilibrium in these two classes of markets.
More details are given in Section~\ref{sect:CES}.
We note that the tatonnement for Leontief Fisher markets that was analysed in~\cite{CCD2013} has an unrealistic constraint on the step sizes;
our analysis removes that constraint, and works for both synchronous and asynchronous tatonnement.
\section{Asynchronous Gradient Descent Model}\label{sect:model}

We consider the following unconstrained optimization problem:
given a convex function $\phi$: $\rr^n\rightarrow \rr$, find its minimal point.
In our model, time, denoted by $t$, is continuous.
The gradient descent process starts at $t=0$ from an initial point $\pc = \left(\pc_1,\pc_2\,\cdots,\pc_n\right)$.
For simplicity, we assume that there are $n$ cores, and $p_j$ is updated by the $j$-th core.\footnote{If
there are fewer cores it suffices to cluster coordinates.}
After each update, the updating core broadcasts it; the other cores receive the message, possibly with a delay.

\noindent\textbf{Notational Convention}~~When there is an update at time $t$ which updates the value of one or more variables,
for each such variable $\square$, we let both $\square^{t-}$ and $\square^t$ denote its value
just before the update,
and $\square^{t+}$ its value right after the update.

We define $p^t \equiv p^{t-}$, the \emph{current point} at time $t$, to comprise the most recently updated values for each coordinate.
However, any particular core may have out-of-date values for one or more coordinates,
but not too much out-of-date, as we specify next.

Let $t_1$ and $t_2$ be the times of successive updates to $p_j$.
Then, at time $t_2$, the $j$-th core will have values for each of the other coordinates that were current at time $t_1$ or later.
In other words, the time taken to communicate an update is no larger than $t_2 - t_1$.
Effectively, this is the constraint on how much parallelism is possible.
Informally speaking, the information which the core holds is at most one ``round'' out of date w.r.t.~its updates.
In fact, it seems likely that we could extend our analysis to allow for any fixed constant number of rounds of datedness,
but as this would entail a proportionate reduction in the step sizes, it does not seem useful.

However, there is no requirement that updates occur at a similar rate, although we imagine that this would be the typical case.
It may be natural in some settings for coordinates to adjust with different frequencies,
e.g.~prices of different goods in a broad enough market.
Accordingly, we define a rather general update rule, as follows.
Each core has the freedom to determine the time at which it updates its coordinate.
To proceed, it will be helpful to define the following rectangular subsets of coordinate values.

\begin{defn}
$\PP{j}{t_1}{t_2}{s_j}$
comprises the rectangular box with $p_j = s_j$ and, for $k\neq j$, spanning the range of values $p_k$
that occur over the time interval $[t_1,t_2]$.
\end{defn}

Let $\tau_j$ be the time at which the last update to $p_j$ occurred,
and let $t$ be the time of the current update to $p_j$.
To update $p_j$, the $j$-th core computes $\nabla_j \phi(\tp)$, where $\tp$ is an arbitrary point in $\PP{j}{\tau_j}{t}{p_j^t}$.
This flexibility allows different coordinates at the $j$-th core to be \emph{variously dated},
under the constraint that they are all no older than time $\tau_j$.
The general form of an update is
$$p_j \leftarrow p_j + F_j(\tp,\nabla_j \phi(\tp),t)\cdot (t-\tau_j),$$
where $F_j$ is a function such that $F_j(\tp,\nabla_j\phi(\tp),t)$ has the same sign as $-\nabla_j\phi(\tp)$.

The term $t-\tau_j$ is somewhat unusual. It is needed because we impose no bound on the frequency of updates.
Without this multiplier, a core, the $k$-th core say, could perform many updates in the time interval $[\tau_j, t]$,
potentially making a cumulatively large update to $p_k$,
which could lead to an unbounded difference between $\nabla_j\phi(\tp)$ and $\nabla_j\phi(p^t)$.
This appears to preclude the usual approaches to a proof of convergence, and even calls convergence into question in general.
If, in fact, $t - \tau_j = \Theta(1)$ always, then this term can be omitted.

Note that the sign of $F_j(\tp,\nabla\phi(\tp),t)$ can be opposite to that of $F_j(p^t,\nabla_j\phi(p^t),t)$;
when this occurs, an update will increase the value of $\phi$, i.e.~we have a bad update!

We do not require any further coordination between the cores.
We just require a minimal amount of communication to ensure that the cores know
an approximation of the current point so that they can compute a useful gradient.
\section{Amortized Analysis}\label{sect:analysis}

\newcommand{\bH}{H}
\newcommand{\Hjk}{H_{jk}}
\newcommand{\Hkj}{H_{kj}}
\newcommand{\tg}{\tilde{g}}
\newcommand{\Dt}{\Delta t}
\newcommand{\Dp}{\Delta p}
\newcommand{\ptauj}{p^{\tau_j}}
\newcommand{\bgam}{\overline{\gamma}}
\newcommand{\tcur}{t^{\bullet}}
\newcommand{\pcur}{p^{\bullet}}
\newcommand{\gjcur}{g_j^\bullet}
\newcommand{\tgjM}{\tg_{j,\max}}
\newcommand{\tgjm}{\tg_{j,\min}}
\newcommand{\gjt}{\gamma_j^t}
\newcommand{\Gjt}{\Gamma_j^t}
\newcommand{\ione}{{i_1}}
\newcommand{\itwo}{{i_2}}
\newcommand{\pj}{p_j}
\newcommand{\pk}{p_k}
\newcommand{\sigj}{\sigma_j}
\newcommand{\sigk}{\sigma_k}
\newcommand{\UH}[5]{\bH_{#1 #2}^{[#3,#4]}\left(#5\right)}

Let $\phi:\rr^n\rightarrow \rr$ be a twice-differentiable convex function.
Our AGD algorithm solves the problem of finding (or approximating) a minimal point of $\phi$, which we denote by $\ps$.
WLOG, we assume that $\phi^* := \phi(\ps) = 0$. We assume that no two updates occur at the same time.\footnote{If
two or more updates do occur at the same time, our analysis remains valid by making infinitesimal perturbations to their update times.}

By default, each core possesses the most up-to-date entry for the coordinate it updates.
However, due to communication delay, it may have outdated entries for coordinates updated by other cores.
Recall that $p^t$ denotes the most up-to-date entries at time $t$;
let $\tp_k^{j,t}$ denote the entry for $p_k$ that the $j$-th core possesses at time $t$.
Note that $\tp^{j,t} \in \PP{j}{\tau_j}{t}{p_j^t}$.

We now consider an update to $p_j$ at time $t$ given by
\begin{equation}\label{eq:update-rule}
p_j' \leftarrow  p_j - \frac{\tg_j(t)}{\gjt} \Dt_j,
\end{equation}
where $\tg_j(t) = \nabla_j\phi(\tp^{j,t})$, $\Dt_j = t - \tau_j$, and
$1/\gjt$ is the \emph{step size}, which will be determined by a rule we specify later.
We assume that $\Dt_j \leq 1$ always, i.e.~two consecutive updates to the same coordinate occur at most one time unit apart.
We note that Rule~\eqref{eq:update-rule} is quite general for it allows both additive and multiplicative updates,
depending on the choice of the $\gjt$. As we shall see, our analysis handles applications of both types.

For any $S\subset \rr^n$, let $\bH_{k\ell}(S) := \max_{p'\in S} \left|\frac{\partial^2 \phi}{\partial p_k \partial p_\ell}(p')\right|$.
We will use the shorthand $\UH{k}{\ell}{t_1}{t_2}{s_\ell}$ for $\bH_{k\ell}\left(\PP{\ell}{t_1}{t_2}{s_\ell}\right)$.
In order to show our convergence results, the $\gjt$ need to be suitably constrained
and the Hessian entries need to be sufficiently bounded. We capture this in our definition of
\emph{controlled} $\gjt$  and $\Hjk$, given right after Theorem~\ref{thm:main} below.

\begin{theorem}\label{thm:main}
Suppose that all updates are made according to update rule \eqref{eq:update-rule}.
Let $\bgam = \max_{j,t} \gamma_j^t$.
If the variables $\gjt$  and $\Hjk$ are controlled, then

\begin{enumerate}
\item[(a)] Suppose the set $\left\{p' \, | \, \phi(p')\leq 2\phi(p^0)\right\}$ is bounded with diameter $B$.
Let $M(B) := \Theta(B^2 \bgam)$.
Then, if $\phi(p^0) \leq M(B)$, $\phi(p^t) = O\left( \frac{M(B)}{t} \right)$;
and otherwise, for $t\leq t' = O\left( \log \frac {\phi(p^0)} {M(B)} \right)$,
$\phi(p^t) = O\left(2^{-\Theta(t)}\phi(p^0)\right)$,
and for $t > t'$, $\phi(p^t) = O\left( \frac{M(B)}{t - t'} \right)$.
\item[(b)] If $\phi$ is strongly convex with parameter $c$,\footnote{i.e.~for any $p_1,p_2$ in its domain,
$\phi(p_2) \geq \phi(p_1) + \nabla \phi(p_1) \cdot (p_2-p_1) + \frac{c}{2} \|p_2 - p_1\|^2$.}
then $\phi(p^t) \leq \left(1-\Theta\left(\frac{c}{\bgam}\right)\right)^t\cdot \phi(p^0)$.
\end{enumerate}
\end{theorem}

\begin{defn}\label{def:control}
The variables $\gjt$ and $\Hjk$ are said to be \emph{controlled} if there
are constants $\alpha \ge 2$, $\eF,\eB > 0$, with $\frac{1}{\alpha} + 2\eB + 2\eF < 1$,
and for each $j$ and time $t$ at which $p_j$ is updated, there are positive numbers $\{\xi_k^t\}_{k\neq j}$,
such that:
\begin{enumerate}
\item[A1.] \emph{(Local Lipschitz bound.)}Let $S_j =\mbox{\emph{Span}}\left\{p_j^{t-},p_j^{t+} \right\}$.
For any $p' \in p^t_{-j} \times S_j$,
$$\phi(p') - \phi(p^t) - \nabla_j\phi(p^t)\cdot (p'_j - p_j^t) \leq \frac {\gjt}{\alpha} (p'_j - p_j^t)^2.$$
\item[A2.] \emph{(Upper bound on $\gjt$.)} For each $j$, there exists a finite positive number $\bgam_j$ such that for all $t$ at which an update to $p_j$ occurs, $\gjt \leq \bgam_j$. We let $\bgam := \max_j \bgam_j$.
\item[A3.] \emph{(Bound on nearby future Hessian entries.)} $ \sum_{k\neq j} \xi_k^t \cdot \UH{j}{k}{t}{\sigk}{p_k^{\tau_k+}} \le \eF \gjt$,
where $\sigk > t$ is the time of the next update to $\pk$;
\item[A4.] \emph{(Bound on recent past Hessian entries.)} $\sum_{k\neq j} \left(\max_{i:k_i=k} \frac{1}{\xi_j^{\beta_i}}\right)\cdot
\UH{k}{j}{\tau_j}{t}{p_j^t} \le \eB \gjt$,
where the index $i$ runs over all updates to coordinate $k$ between times $\tau_j$ and $t$, and $\beta_i$ is the time at which each such update occurs
(this notation is defined precisely in Lemma \ref{lem:error-of-gradient}).
\end{enumerate}
\end{defn}

If the updates used fully up-to-date gradients,
i.e.\ if $\Dp_j = -\frac{\nabla_j\phi(p^t)}{\gjt} \Dt_j$,
rearranging Condition A1 would give the following lower bound on the progress (cf.~Lemma~\ref{lem:progress-of-phi} below):
$$\phi(p^{t-}) - \phi(p^{t+}) \ge \sum_j \frac {1}{\gjt} (\nabla_j\phi(p^t))^2  \Dt_j -  \frac {1}{\alpha\gjt} (\nabla_j\phi(p^t))^2 \Dt_j^2 \ge \sum_j \left( 1 - \frac{1}{\alpha} \right) \frac{(\nabla_j\phi(p^t))^2 \Dt_j}{\gjt}.$$
The remaining conditions are present to cope with the lack of synchrony.
Conditions A3 and A4 ensure that the ``errors' in the gradients we use for the updates
are not too large cumulatively. Basically, they will reduce the multiplier in the progress
from $(1 - \frac {1}{\alpha})$ to $(1 - \frac {1}{\alpha} -2\eF - 2\eB)$.
Recall that the lack of synchrony may result in bad updates.
To hide the resulting temporary lack of progress and to show continued long-term progress,
we use an amortized analysis which employs the following potential function.
\begin{equation}\label{eq:PF}
\Phi(p^t,t,\tau) = \phi(p^t) - c_1 \sum_j \int_{\tau_j}^t \frac{(g_j(t'))^2}{\bgam_j}\,dt' + \sum_j \sum_i \xi^{\beta_i}_j \cdot \UH{k_i}{j}{\beta_i}{\sigj}{p_j^{\tau_j+}}\frac{\left(\Dp_{k_i}\right)^2}{\Dt_{k_i}}\left[2-c_2 (t - \beta_i)\right],
\end{equation}
where $g_j(t') := \nabla_j \phi\left(p^{t'}\right)$ and $\sigj > \tau_j$ is the time of the next update to $\pj$;
for each $j$, the index $i$ runs over all updates, between times $\tau_j$ and $t$,
to coordinates other than $j$;
$c_1$ and $c_2$ are positive constants whose values we will determine later.
$\left\{\xi^{\beta_i}_j\right\}$ are the positive numbers in Conditions A3 and A4;
note that these variables are indexed by $i$ but not by the update coordinate $k_i$,
so for any $j$, $\xi^{\beta_\ione}_j$ may be different from $\xi^{\beta_\itwo}_j$, even if $k_\ione = k_\itwo$.

The integral in the above potential function reflects the ideal progress were there a continuous synchronized
updating of the prices, and the additional terms are present to account for the attenuation of progress
due to asynchrony.

Our method of analysis is to show that $\frac{d\Phi}{dt} \leq - \beta_1 \Phi^2$ for a suitable
constant $\beta_1>0$ whenever there is no price update, and that $\Phi$ only decreases when there is a
price update; this then yields Theorem \ref{thm:main}(a). Theorem \ref{thm:main}(b) follows from a stronger bound on the
derivative, namely that $\frac{d\Phi}{dt} \le - \beta_2 \Phi$, where $\beta_2>0$.
This general approach for asynchrony analysis was used previously by Cheung et al.~\cite{CCR2012} for a result in the style of (b), but for a quite different potential function.

\bigskip

It is straightforward to show that when there is no update,
\begin{equation}\label{eq:cont-progress}
\frac{d\Phi}{dt} = -c_1 \sum_j \frac{(g_j(t))^2}{\bgam_j} - c_2 \sum_j \sum_i \xi^{\beta_i}_j \cdot \UH{k_i}{j}{\beta_i}{\sigj}{p_j^{\tau_j+}} \frac{\left(\Dp_{k_i}\right)^2}{\Dt_{k_i}}.
\end{equation}

Lemma~\ref{lem:progress-of-phi} below bounds the change to $\phi$ when there is an update.
Lemma~\ref{lem:error-of-gradient} states some useful bounds on the maximum change that can occur to the gradient between
two updates to the same coordinate.
Lemma~\ref{lem:Phi-at-update} below bounds the change to $\Phi$ when there is an update.

\newcommand{\PPP}[1]{\left(\PP{j}{#1}{t}{p_j^t}\right)}
\newcommand{\bareta}{\bar{\eta}}

\begin{lemma}\label{lem:progress-of-phi}
Suppose there is an update to $p_j$ at time $t$ according to rule \eqref{eq:update-rule}, with $\gjt$ satisfying Condition A1.
Let $\phi^-$ and $\phi^+$ denote, respectively, the convex function values just before and just after the update.
Let $g_j := \nabla_j \phi(p^t)$ and $\tg_j \equiv \tg_j(t)$.
Let $\Dp_j$ be the change to $p_j$ made by the update, i.e.~$\Dp_j := - \frac{\tg_j(t)}{\gjt} \Dt_j$.
Then
$$\phi^- - \phi^+ \geq \left(1-\frac{1}{\alpha}\right)\frac{\gjt (\Dp_j)^2}{\Dt_j} - |g_j - \tg_j| \cdot |\Dp_j|.$$
\end{lemma}

\begin{lemma}\label{lem:error-of-gradient}
Suppose that between times $\tau_j$ and $t$, there are updates to the sequence of coordinates $k_1,k_2,\cdots,k_m$,
which may include repetitions,
but include no update to coordinate $j$.
Let $\beta_1,\beta_2,\cdots,\beta_m$ denote the times at which these updates occur.
Let $\tgjM$ and $\tgjm$ denote, respectively, the maximum and minimum values of $\nabla_j(p')$, where $p'\in \PP{j}{\tau_j}{t}{p_j^t}$.
For any positive numbers $\{\eta_i\}_{i=1\cdots m}$, for each $k\neq j$, let $\bareta_k := \min_{i:k_i = k} \eta_i$.
Then for any real number $\mu$,
\begin{equation}\label{eq:error-of-gradient-2}
|\mu| \cdot \left(\tgjM - \tgjm\right) \leq 2 \mu^2\sum_{k\neq j} \frac{1}{\bareta_k} \UH{k}{j}{\tau_j}{t}{p_j^t}
 + \sum_{i=1}^m \eta_i \cdot \UH{k_i}{j}{\beta_i}{t}{p_j^t} \frac{(\Dp_{k_i})^2}{\Dt_{k_i}}
\end{equation}
and
\begin{equation}\label{eq:error-of-gradient-3}
\left(\tgjM - \tgjm\right)^2 \leq 8 \left(\sum_{i=1}^m \eta_i \cdot \UH{k_i}{j}{\beta_i}{t}{p_j^t}\frac{(\Dp_{k_i})^2}{\Dt_{k_i}}\right)
\left(\sum_{k\neq j} \frac{1}{\bareta_k} \UH{k}{j}{\tau_j}{t}{p_j^t}\right).
\end{equation}
\end{lemma}

\begin{lemma}\label{lem:Phi-at-update}
Suppose that there is an update to $p_j$ at time $t$. Suppose that $\gjt$ is chosen so that Conditions A1, A3 and A4 hold.
Let $\Phi^-$ and $\Phi^+$, respectively, denote the values of $\Phi$ just before and just after the update. Then
\begin{align*}
\Phi^- - \Phi^+ & \geq \left(1-\frac{1}{\alpha} - 2\eB - c_1(1+4\eB) - 2\eF\right)\frac{\gjt (\Dp_j)^2}{\Dt_j}\\
&\qquad + \left(1 - c_2 - c_1(2+8\eB)\right) \sum_{i=1}^m \xi_j^{\beta_i} \cdot \UH{k_i}{j}{\beta_i}{t}{p_j^t} \frac{(\Dp_{k_i})^2}{\Dt_{k_i}}.
\end{align*}
\end{lemma}

\begin{pf}
By Lemma~\ref{lem:progress-of-phi} and the fact $(t-\beta_i)\leq (t-\tau_j) \leq 1$,
\begin{align}
\Phi^- - \Phi^+ &= \phi^- - \phi^+ - c_1 \int_{\tau_j}^t \frac{(g_j(t'))^2}{\bgam_j}\,dt' + \sum_i \xi^{\beta_i}_j \cdot \UH{k_i}{j}{\beta_i}{t}{p_j^{\tau_j+}}
\frac{(\Dp_{k_i})^2}{\Dt_{k_i}}\left[2-c_2(t-\beta_i)\right] \nonumber\\
&\qquad\qquad - 2 \sum_{k\neq j} \xi^t_k \cdot \UH{j}{k}{t}{\sigk}{p_k^{\tau_k+}}\frac{(\Dp_j)^2}{\Dt_j} \nonumber\\
&~\geq \left(1-\frac{1}{\alpha}\right)\frac{\gjt (\Dp_j)^2}{\Dt_j}
- \underbrace{|g_j - \tg_j| \cdot |\Dp_j|}_{E_1}
- \underbrace{c_1 \int_{\tau_j}^t \frac{(g_j(t'))^2}{\bgam_j}\,dt'}_{E_2} \nonumber\\
&\qquad\qquad + (2-c_2)\sum_i \xi^{\beta_i}_j \cdot \UH{k_i}{j}{\beta_i}{t}{p_j^{\tau_j+}}\frac{(\Dp_{k_i})^2}{\Dt_{k_i}}
- \underbrace{2 \sum_{k\neq j} \xi^t_k \cdot \UH{j}{k}{t}{\sigk}{p_k^{\tau_k+}}\frac{(\Dp_j)^2}{\Dt_j}}_{E_3}.\label{eq:change_of_Phi_inter}
\end{align}

We bound $E_1,E_2$ and $E_3$ below. We will be applying \eqref{eq:error-of-gradient-2} and \eqref{eq:error-of-gradient-3} with $\eta_i = \xi^{\beta_i}_j$.
Let
$$V_1 := \sum_{k\neq j} \frac{1}{\min_{i:k_i=k} \xi^{\beta_i}_j} \UH{k}{j}{\tau_j}{t}{p_j^t}\qquad\mbox{and}\qquad
V_2 := \sum_{i=1}^m \xi^{\beta_i}_j\cdot \UH{k_i}{j}{\beta_i}{t}{p_j^t} \frac{(\Dp_{k_i})^2}{\Dt_{k_i}}.$$
Note that by Condition A4, $V_1 \leq \eB \gjt$. By \eqref{eq:error-of-gradient-2},
$E_1 \leq 2(\Dp_j)^2 V_1 + V_2 \leq 2 \eB \gjt (\Dp_j)^2 + V_2.$

To bound $E_2$, first note that for any $t'\in (\tau_j,t]$, $p^{t'}\in \PP{j}{\tau_j}{t}{p_j^t}$. Then
\begin{align}
& \frac{(g_j(t'))^2}{\bgam_j} - \frac{(\tg_j)^2}{\bgam_j} = \frac{(g_j(t') - \tg_j)^2}{\bgam_j} - \frac{2 \tg_j}{\bgam_j} (\tg_j - g_j(t'))\nonumber\\
 & \leq~ \frac{(g_j(t') - \tg_j)^2}{\bgam_j} + 2 \left|\frac{\tg_j}{\bgam_j}\right| \cdot \left|\tg_j - g_j(t')\right|
 ~\leq~ \frac{8}{\bgam_j} V_2 V_1 + \frac{4(\tg_j)^2}{(\bgam_j)^2} V_1 + 2V_2\comm{by Eqns.~\eqref{eq:error-of-gradient-3} and \eqref{eq:error-of-gradient-2}}\nonumber\\
 & \leq~ \frac{8\eB \gjt}{\bgam_j} V_2 + \frac{4\eB \gjt (\tg_j)^2}{(\bgam_j)^2} + 2V_2
~\leq~\frac{4\eB (\tg_j)^2}{\bgam_j} + (2+8\eB) V_2 \comm{by Condition A2}\label{eq:bound-integral-1}
\end{align}
Hence $\frac{(g_j(t'))^2}{\bgam_j} \leq (1+4\eB) \frac{(\tg_j)^2}{\bgam_j} + (2+8\eB) V_2$,
and then as $\Dt_j\leq 1$,
\begin{equation}\label{eq:bound-integral-2}
E_2 ~\leq~  {c_1}\int_{\tau_j}^t \frac {(g_j(t'))^2} {\bgam_j}  dt'
~\leq~ c_1(1+4\eB)\frac{(\tg_j)^2 \Dt_j}{\bgam_j} + c_1(2+8\eB) V_2
\end{equation}

Finally, by Condition A3, $E_3 \leq 2 \eF \gjt \frac{(\Dp_j)^2}{\Dt_j}$.

\smallskip
\noindent Combining the above bounds on $E_1,E_2,E_3$ yields
\begin{align*}
\Phi^- - \Phi^+ &\geq \left(1-\frac{1}{\alpha}\right)\frac{\gjt (\Dp_j)^2}{\Dt_j} - \left[2 \eB\gjt (\Dp_j)^2 + V_2\right]
 - \left[c_1(1+4\eB)\frac{(\tg_j)^2 \Dt_j}{\bgam_j} + c_1(2+8\eB) V_2\right]\\
&\qquad + (2-c_2)V_2 - 2 \eF \gjt \frac{(\Dp_j)^2}{\Dt_j}.
\end{align*}
As $\Dp_j = - \frac{\tg_j(t)}{\gjt} \Dt_j$ and $\Dt_j \leq 1$, the result follows.
\end{pf}

\begin{lemma}\label{lem:relating-Phi-phi}
If $2 -c_2 \geq c_1(2+8\eB)$, then $\Phi(p^t,t,\tau) \geq \left[1-2c_1(1+4\eB)\right]\phi(p^t)$.
\end{lemma}

\begin{pfof}{Theorem \ref{thm:main}(a)}
Choose
$c_1 = (1+4\eB)^{-1}\cdot \min\left\{ 1 - \frac{1}{\alpha} - 2\eB - 2\eF, \frac{1}{4} \right\}$ and $c_2 = 1 - c_1(2+8\eB)$.
Then the following hold:
(i) $c_1,c_2 > 0$;
(ii) $1 - \frac{1}{\alpha} - 2\eB - 2\eF - c_1(1+4\eB) \geq 0$;
(iii) $1 - c_2 - c_1(2+8\eB) = 0$;
(iv) $2-c_2 \geq c_1(2+8\eB)$;
(v) $c_1(1+4\eB) \leq \frac{1}{4}$.

By (ii), (iii) and Lemma \ref{lem:Phi-at-update}, $\Phi$ does not increase at any update.

By (iv), (v) and Lemma \ref{lem:relating-Phi-phi}, $\Phi(p^t,t,\tau) \geq \frac{\phi(p^t)}{2}$.
Thus, $\forall t\geq 0$, $\phi(p^t) \leq 2\Phi(p^t,t,\tau)\leq 2\Phi(p^0,0,\vec{0}) = 2\phi(p^0)$,
i.e.~$\{p^t\}_{t\geq 0}$ is contained in the set $\{p'\,|\,\phi(p')\leq 2\phi(p^0)\}$,
which, by assumption, has diameter at most $B$.

Note that at any time $t$, by the convexity of $\phi$, $\phi(p^t) + \sum_j g_j(t) \cdot (p_j^* - p_j^t) \leq \phi^* = 0$ and hence
$$\sum_j |g_j(t)| \cdot |p_j^t - p_j^*| \geq \sum_j g_j(t) \cdot (p_j^t - p_j^*) \geq \phi(p^t) \geq 0.$$
By the Cauchy-Schwarz inequality,

$$\phi(p^t) \leq \sum_j |g_j(t)| \cdot |p_j^t - p_j^*| \leq \sqrt{\left(\sum_j (g_j(t))^2\right)\left(\sum_j (p_j^t - p_j^*)^2\right)} \leq B \sqrt{\sum_j (g_j(t))^2}.$$
Then
$$\sum_j \frac{(g_j(t))^2}{\bgam_j} \geq \frac{1}{\bgam}\sum_j (g_j(t))^2 \geq \frac{1}{\bgam}\left(\frac{\phi(p^t)}{B}\right)^2 = \frac{1}{B^2 \bgam} \phi(p^t)^2.$$
By \eqref{eq:cont-progress},
$$\frac{d\Phi}{dt} \leq -\frac{c_1}{B^2\bgam}\cdot \phi(p^t)^2 - c_2 \sum_j \sum_i \xi^{\beta_i}_j \cdot \UH{k_i}{j}{\beta_i}{\sigj}{p_j^{\tau_j+}}
\frac{(\Dp_{k_i})^2}{\Dt_{k_i}}.$$

By \eqref{eq:PF}, $\Phi(p^t,t,\tau) \leq \phi(p^t) + 2 \sum_j \sum_i \xi^{\beta_i}_j \cdot \UH{k_i}{j}{\beta_i}{\sigj}{p_j^{\tau_j+}} \frac{(\Dp_{k_i})^2}{\Dt_{k_i}}$. Let $X_1 := \phi(p^t)$ and $X_2 := \sum_j \sum_i \xi^{\beta_i}_j \cdot \UH{k_i}{j}{\beta_i}{\sigj}{p_j^{\tau_j+}} \frac{(\Dp_{k_i})^2}{\Dt_{k_i}}$. Then $\Phi \leq X_1 + 2 X_2$ and $\frac{d\Phi}{dt} \leq -\frac{c_1}{B^2 \bgam} (X_1)^2 - c_2 X_2$. Let $M(B) := \Theta(B^2 \bgam)$.
As $\phi(p^t) \le 2\Phi(t)$, this guarantees that
if $\phi(p^0) = \Phi(p^0) \leq M(B)$, then $\phi(p^t) = O\left( \frac{M(B)}{t} \right)$;
and otherwise, for $t\le  t' = O\left( \log \frac {\phi(p^0)} {M(B)} \right)$,
$\phi(p^t) = O\left(2^{-\Theta(t)}\phi(p^0)\right)$,
and for $t > t'$, $\phi(p^t) = O\left( \frac{M(B)}{t - t'} \right)$.
\end{pfof}

\begin{pfof}{Theorem \ref{thm:main}(b)}
If $\phi$ is strongly convex with parameter $c$, then, by definition,
\begin{align*}
0 = \phi^* & \geq \phi(p^t) + \sum_j g_j(t) \cdot (p^*_j - p_j^t) + \frac{c}{2} \sum_j (p^*_j - p_j^t)^2 \\
&\geq \phi(p^t) + \min_{p'}\left\{\sum_j g_j(t) \cdot (p'_j - p_j^t) + \frac{c}{2}(p'_j - p_j^t)^2\right\}.
\end{align*}
Computing the minimum point of the quadratic polynomial in $(p'_j - p_j^t)$ yields $0 \geq \phi(p^t) - \sum_j \frac{(g_j(t))^2}{2c}$.
Then
$$\sum_j \frac{(g_j(t))^2}{\bgam_j} \geq \frac{1}{\bgam}\sum_j (g_j(t))^2 \geq \frac{2c}{\bgam}\phi(p^t).$$
As in Case (a), $\Phi \leq X_1 + 2 X_2$; and by \eqref{eq:cont-progress}, $\frac{d\Phi}{dt} \leq -\frac{2cc_1}{\bgam} X_1 - c_2 X_2$.
This guarantees that $2\phi(p^t) \le \Phi(t) \leq (1-\delta(c))^t \phi(p^0)$, where $\delta(c) = \min \{ \frac {c c_1} {\bgam}, \frac {c_2}{4} \}$.
\end{pfof}
\section{Solving Matrix Systems}\label{sect:linear-systems}

\newcommand{\transpose}{^{\tiny \mbox{T}}}

For any symmetric and positive definite (SPD) matrix $A\in \rr^{n\times n}$ and $b,p\in \rr^n$,
let $f_{A,b} (p) = \frac{1}{2} p\transpose Ap - p\transpose b$.
It is well known that $f_{A,b}(p)$ is a strictly convex function of $p$, and $\nabla f_{A,b}(p) = Ap - b$.
Therefore, finding the minimum point of $f_{A,b}(p)$ is equivalent to solving the linear system $Ap = b$,
and hence one can solve the linear system by performing gradient descent on $f_{A,b}(p)$.

The Hessian of $f_{A,b}(p)$ is $\nabla^2 f_{A,b}(p) = A$, a constant matrix.
This allows a simple rule to determine a \emph{constant} step size for each coordinate.
By taking all the $\xi$ values to be $1$, to apply Theorem~\ref{thm:main},
it suffices to have $\gjt = \gamma_j$ satisfy
$\gamma_j \geq \frac{A_{jj}}{2} \alpha$ (for A1), $\frac{4}{\gamma_j} \sum_{k\neq j} |A_{jk}| < 1 - \frac{1}{\alpha}$ (combining A3, A4 and the bound
$\frac {1}{\alpha} + 2\eF + 2\eB < 1$),
and $\alpha\geq 2$.
These imply it suffices that the step size, $1/\gamma_j$, be less than
$\left[\max\left\{\frac{A_{jj} + 8\sum_{k\neq j} |A_{kj}|}{2},A_{jj}\right\}\right]^{-1}$.

Another application is given by the following class of optimization problems (see Nesterov~\cite{Nesterov2012}): minimizing
$F(p) := \sum_{i=1}^n f_i(p_i) + \frac{1}{2} \|Ap - b\|^2$,
where the $f_i$ are convex differentiable univariate functions,
$A\in \rr^{r\times n}$ is an $r\times n$ real matrix and $b\in \rr^r$.
The Hessian of $F$ at $p$ is $A\transpose A+D$,
where $D$ is the diagonal matrix with $D_{jj} = f_j ''(p_j)$.
If $f_j''(p)$ is bounded by $L_j$, again, it suffices to have $\gjt = \gamma_j$ satisfy
$\gamma_j \geq \frac{(A\transpose A)_{jj} + L_j}{2} \alpha$,
$\frac{4}{\gamma_j} \sum_{k\neq j} |(A\transpose A)_{jk}| < 1 - \frac{1}{\alpha}$, and $\alpha\geq 2$.
These imply it suffices that the step size, $1/\gamma_j$, be less than $\left[\max\left\{\frac{(A\transpose A)_{jj} + L_j + 8\sum_{k\neq j} |(A\transpose A)_{kj}|}{2},(A\transpose A)_{jj} + L_j\right\}\right]^{-1}$.

Next, we discuss how $\nabla_j F(p)$ is computed by the $j$-th core.
Let $G(p) = Ap-b$ and let $A_j$ denote the $j$-th column of the matrix $A$. Then $\nabla_j F(p) = f_j'(p_j) + (A_j)\transpose G(p)$.
$f_j'(p_j)$ is recomputed only when $p_j$ changes.
For any $k$, when $p_k$ is changed by $\Dp_k$, $G(p+\Dp_k) - G(p) = \Dp_k A_k$, and hence
$(A_j)\transpose G(p)$ changes by $\Dp_k (A_j)\transpose A_k$.
Note that $(A_j)\transpose A_k$ is a constant and hence can be pre-calculated,
so the above equation provides a quick way to update $\nabla_j F(p)$ once the $j$-th core receives the message with $\Dp_k$.

Recall that our AGD algorithm allows different coordinate values to be variously dated,
under the constraint that they are all no older than the time of the last update.
It is natural to aim to have essentially the same frequency of update for each coordinate.
Accordingly, at the $i$-th round of updates, each core can simply ensure
it has received the update for the previous round from every other core.
The update messages might arrive at different times, 
but the $j$-th core needs not wait until it collects all such messages.
It can simply compute the changes to $\nabla_j F(p)$ incrementally
as it receives updates $\Delta p_k$ to $p_k$.
This avoids the need for any explicit synchronization.
\section{Tatonnement in Fisher Markets}\label{sect:CES}

A \emph{Fisher market} comprises a set of $n$ goods and two sets of agents, sellers and buyers.
The sellers bring the goods to market and the buyers bring money with which to buy the goods.
The trade is driven by a collection of non-negative prices $\{p_j\}_{j=1\cdots n}$, one price per good.
WLOG, we assume that each seller brings one distinct good to the market, and she is the price-setter for this good.
By normalization, we may assume that each seller brings one unit of her good to the market.

Each buyer $i$ starts with $e_i$ money, and has a utility function
$u_i(x_{i1},x_{i2}, \cdots ,x_{in})$ expressing her preferences:
if she prefers bundle $\{x^a_{ij}\}_{j=1\cdots n}$ to bundle $\{x^b_{ij}\}_{j=1\cdots n}$, then
$u_i ( \{x^a_{ij}\}_{j=1\cdots n} ) > u_i ( \{x^b_{ij}\}_{j=1\cdots n} )$.
At any given prices $\{p_j\}_{j=1\cdots n}$, each buyer $i$ seeks to purchase a maximum utility bundle of goods costing at most $e_i$.
The \emph{demand} for good $j$, denoted by $x_j$, is the total quantity of the good sought by all buyers.
The \emph{supply} of good $j$ is the quantity of good $j$ its seller brings to the market, which we have assumed to be $1$.
The \emph{excess demand} for good $j$, denoted by $z_j$, is the demand for the good minus its supply, i.e.~$z_j = x_j - 1$.
Prices $\{\ps_j\}_{j=1\cdots n}$ are said to form a \emph{market equilibrium} if,
for any good $j$ with $\ps_j>0$, $z_j = 0$,
and for any good $j$ with $\ps_j=0$, $z_j \leq 0$.

The following two classes of utility functions are commonly used in market models.
The first class is the Constant Elasticity of Substitution (CES) utility function:
$$u_i\left(x_{i1},x_{i2},\cdots,x_{in}\right) = \left(a_{i1} (x_{i1})^{\rho_i} + a_{i2} (x_{i2})^{\rho_i}
+ \cdots + a_{in} (x_{in})^{\rho_i}\right)^{1/\rho_i},$$
where $\rho_i\leq 1$ and $\forall j$, $a_{ij}\geq 0$.
$\theta_i := \rho_i / (\rho_i - 1)$ is a parameter which will be used in the analysis.
In this paper we focus on the cases $\rho_i\leq 0$,
in which goods are complements and hence the utility function is called a complementary-CES utility function.
It is easy to extend our analysis to the cases $\rho_i\geq 0$, which had been analysed in \cite{CF2008,CFR2010}.
The second class is the Leontief utility function:
$$u_i\left(x_{i1},x_{i2},\cdots,x_{in}\right) = \min_{j\in S} \left\{b_{ij} x_{ij}\right\},$$
where $S$ is a non-empty subset of the goods in the market, and $\forall j\in S$, $b_{ij} > 0$.

Cheung, Cole and Devanur \cite{CCD2013} showed that
tatonnement is equivalent to gradient descent on a convex function $\phi$ for Fisher markets
with buyers having complementary-CES or Leontief utility functions (defined in the appendix).
To be specific, $\nabla_j\phi(p) = -z_j(p)$, and the convex function $\phi$ is
$\phi(p) = \sum_j p_j + \sum_i \hat{u}_i(p)$,
where $\hat{u}_i(p)$ is the optimal utility that buyer $i$ can attain at prices $p$.
The corresponding update rule is
\begin{equation}\label{eq:async-tat-CES-rule}
p_j' = p_j\cdot \left(1 + \lambda \cdot \min\{\tilde{z}_j,1\}\cdot (t-\tau_j)\right),
\end{equation}
where $\tilde{z}_j$ is a value between the minimum and maximum excess demands during the time interval $(\tau_j, t]$,
and $\lambda > 0$ is a suitable constant. As the update rule is multiplicative, we assume that the initial prices are positive.

Note that $\gjt =\frac{\max\{1,\tilde{z}_j\}}{\lambda p_j}$.
As we will see, it suffices that $\lambda \le \frac{1}{23.46}$.
In comparison, in the synchronous version, $\gjt \geq \frac{6 \max\{1,z_j^t\}}{p_j}$,
so the step sizes of the asynchronous tatonnement are a constant fraction of those used in its synchronous counterpart.

\begin{theorem}
\label{thm:CES-Fisher-cnvge}
For  $\lambda \le \frac{1}{23.46}$, asynchronous tatonnement price updates using rule \eqref{eq:async-tat-CES-rule} converge toward
the market equilibrium in any complementary-CES or Leontief Fisher market.
\end{theorem}

In a Fisher market with buyers having complementary-CES utility functions,
Properties 1 and 2 below are well-known.  Property 3
was proved in~\cite{CCD2013} and implies that Condition A1 holds when $\alpha = 6$ and $\gjt \geq 9.5 x_j(p^t) / p_j^t$.
\begin{enumerate}
\item Let $x_{i\ell}(p)$ denote the buyer $i$'s demand for good $\ell$ at prices $p$. Then for $k\neq j$,
$$\left|\frac{\partial^2 \phi}{\partial p_j \partial p_k}\right| = \sum_i \frac{\theta_i x_{ij}(p) x_{ik}(p)}{e_i} \leq \sum_i \frac{x_{ij}(p) x_{ik}(p)}{e_i}.$$
\item Given positive prices $p$, for any $0 < r_1 < r_2$,
let $p'$ be prices such that for all $j$, $r_1 p_j \leq p'_j \leq r_2 p_j$.
Then for all $j$, $\frac{1}{r_2} x_j(p) \leq x_j(p') \leq \frac{1}{r_1} x_j(p)$.
\item If $\frac{\Dp_j}{p_j} \leq 1/6$, then $\phi(p+\Dp) - \phi(p) - \nabla_j\phi(p) \cdot \Dp_j \leq \frac{1.5 x_j}{p_j} (\Dp_j)^2$.
\end{enumerate}

We outline the analysis for the complementary-CES case.
As $\lambda \le \frac{1}{23.46}$, within one unit of time,
each price can vary by a factor between $(9/10)^2 = 81/100$ and $(11/10)^2 = 121/100$.\footnote{These bounds are loose, but they suffice for our purpose.}
Hence, within one unit of time, the demand can vary by a factor between $100/121$ and $100/81$.

For each update to $p_j$ at time $t$, we choose $\xi_k^t := p_k^t / p_j^t$.
Then the following lemma bounds the sums in Conditions A3 and A4.

\begin{lemma}
\label{lem:ces-bounds}
\emph{(a)} $\sum_{k\neq j} \xi_k^t \cdot \UH{j}{k}{t}{\sigma_j}{p_k^{\tau_k+}} \leq \frac{1.53 x_j(p^t)}{p_j^t}$;\\
\emph{(b)} $\sum_{k\neq j} \left(\max_{q:k_q=k} \frac{1}{\xi_j^{\beta_q}}\right)\cdot \UH{k}{j}{\tau_j}{t}{p_j^t} \leq \frac{1.89 x_j(p^t)}{p_j^t}$.
\end{lemma}

\begin{pf}
\begin{align*}
& \sum_{k\neq j} \xi_k^t \cdot \UH{j}{k}{t}{t+1}{p_k^{\tau_k+}} ~=~ \sum_{k\neq j} \frac{p_k^t}{p_j^t}\cdot \max_{p'\in \PP{k}{t}{t+1}{p_k^{\tau_k+}}}\left|\frac{\partial^2 \phi}{\partial p_j \partial p_k}\right|\\
&~~~~~\leq~ \frac{1}{p_j^t}\sum_{k\neq j} p_k^t \cdot  \max_{p'\in \PP{j}{t}{t+1}{p_k^{\tau_k+}}}\sum_i \frac{x_{ij}(p') x_{ik}(p')}{e_i}
~\leq~ \frac{1}{p_j^t}\sum_{k\neq j} p_k^t \sum_i \frac{\left(\frac{100}{81} x_{ij}(p^t)\right) \left(\frac{100}{81} x_{ik}(p^t)\right)}{e_i}\\
&~~~~~\leq~ \frac{1.53}{p_j^t}\sum_i x_{ij}(p^t) \sum_{k\neq j} \frac{p_k^t x_{ik}(p^t)}{e_i}
~\leq~ \frac{1.53}{p_j^t}\sum_i x_{ij}(p^t) ~=~ \frac{1.53 x_j(p^t)}{p_j^t}.
\end{align*}
And
\begin{align*}
&\sum_{k\neq j} \left(\max_{q:k_q=k} \frac{1}{\xi_j^{\beta_q}}\right)\cdot \UH{k}{j}{\tau_j}{t}{p_j^t} ~=~ \sum_{k\neq j} \frac{\max_{q:k_q=k} p_k^{\beta_q}}{p_j^t} \cdot \max_{p'\in \PP{j}{\tau_j}{t}{p_j^t}}\left|\frac{\partial^2 \phi}{\partial p_j \partial p_k}\right|\\
&~~~~~\leq~ \frac{1}{p_j^t}\sum_{k\neq j} \left(\frac{100}{81} p_k^t\right) \sum_i \frac{\left(\frac{100}{81} x_{ij}(p^t)\right) \left(\frac{100}{81} x_{ik}(p^t)\right)}{e_i}\\
&~~~~~\leq~ \frac{1.89}{p_j^t} \sum_i x_{ij}(p^t) \sum_{k\neq j} \frac{p_k^t x_{ik}(p^t)}{e_i}
~\leq~ \frac{1.89}{p_j^t} \sum_i x_{ij}(p^t) ~=~ \frac{1.89 x_j(p^t)}{p_j^t}.
\end{align*}
\end{pf}

\begin{pfof}{Theorem~\ref{thm:CES-Fisher-cnvge} for the CES case}
By Property 3, Condition A1 is satisfied by setting $\gjt \geq \frac{9.5 x_j(p^t)}{p_j^t}$ and $\alpha=6$.
By Lemma \ref{lem:ces-bounds}, Conditions A3 and A4 are satisfied by setting $\eF = 1/6$ and $\eB = 1/5$,
and $1 - \frac {1} {\alpha} - 2\eF -2 \eB = \frac {1}{10} > 0$.

As discussed in~\cite{CF2008}, the seller might know only $\tilde{x}_j$ but not $x_j$.
As $\tilde{x}_j \geq \frac{81}{100} x_j$,
it would be more natural to use $\gjt \geq \frac{11.73 \tilde{x}_j}{p_j}$,
or the even weaker (but still more natural) $\gjt \geq \frac{23.46 \max\{1,\tilde{z}_j\}}{p_j}$,
which yields update rule \eqref{eq:async-tat-CES-rule}.

\cite{CCD2013} proved that prices in tatonnement cannot get arbitrarily close to zero
and hence demands cannot increase indefinitely,
so $\bgam_j$, as defined in Condition A2, is finite.
\cite{CCD2013} also showed that $\phi$ is strongly convex.
The result follows from Theorem \ref{thm:main}(b).
\end{pfof}

\paragraph{Ongoing Complementary-CES Fisher Markets}
Cole and Fleischer's Ongoing market model \cite{CF2008}
incorporates asynchronous tatonnement and warehouses to form a self-contained dynamic market model.
The price update rule is designed to achieve two goals simultaneously:
convergence toward the market equilibrium and warehouse ``balance''.
As in~\cite{CCR2012},
we modify the price update rule \eqref{eq:async-tat-CES-rule} to achieve both targets.
Analysing its convergence entails the design of a significantly more involved potential function;
the details are given in the appendix.

\paragraph{Leontief Fisher Markets}
It is well-known that Leontief utility functions can be considered as the ``limit'' of CES utility functions as $\rho\ra -\infty$.
Our analysis for CES Fisher markets can be reused, with no modification needed,
to show that in any Leontief Fisher market, $\Phi(p^t,t,\tau)$ decreases with $t$.
However, as an equilibrium price in a Leontief Fisher market can be zero, it is unavoidable that
the chosen step size $\gjt$ may tend to infinity (as $\gjt = \Omega(1/p_j)$),
violating Condition A2; thus Theorem \ref{thm:main} cannot be applied directly.

On top of the result that $\Phi(p^t,t,\tau)$ decreases with $t$, we provide additional arguments to
show that tatonnement with update rule \eqref{eq:async-tat-CES-rule} still converges toward the market equilibrium in Leontief Fisher markets.
The proof is given in the appendix.
However, this result does not provide a bound on the rate of convergence, which appears to preclude incorporating
warehouses into the analysis.

\paragraph{Further Discussion of Asynchronous Dynamics}
Computer science has long been concerned with the organization and manipulation of information in the form of
well-defined problems with a clear intended outcome.
But in the last 15 years, computer science has gained a new dimension,
in which outcomes are predicted or described, rather than designed.
Examples include bird flocking \cite{Chazelle2009}, influence systems \cite{Chazelle2012},
spread of information memes across the Internet \cite{LBK2009} and market economies~\cite{CF2008}.
Many of these problems fall into the broad category of analysing dynamic systems.
Dynamic systems are a staple of the physical sciences;
often the dynamics are captured via a neat, deterministic set of rules
(e.g.~Newton's law of motion, Maxwell's equations for electrodynamics).
The modeling of dynamic systems with intelligent agents presents new challenges because agent behavior may not
be wholly consistent or systematic.
One issue that has received little attention is the timing of agents' actions.
Typically, a fixed schedule has been assumed (e.g.~synchronous or round robin),
perhaps because it was more readily analysed.

This work provides a second demonstration (the first demonstration is in \cite{CFR2010,CCR2012})
and further development of a method for analysing asynchronous dynamics,
here for dynamics which are equivalent to gradient descent.
This methodology may be of wider interest.

\newpage
\bibliographystyle{plain}
\bibliography{agd_bib}

\newpage
\appendix
\def\pjq{p_{j,q}}
\def\tzjq{\tilde z_{j,q}}
\def\D{\Delta}

\section{Missing Proofs in Section \ref{sect:analysis}}

\begin{pfof}{Lemma \ref{lem:progress-of-phi}}
By Condition A1, $\phi^+ - \phi^- - g_j \Dp_j \leq \frac{\gamma_j^t}{\alpha} (\Dp_j)^2$. Then
\begin{align*}
\phi^- - \phi^+ ~&~\geq -[\tg_j + (g_j - \tg_j)] \Dp_j - \frac{\gamma_j^t} {\alpha} (\Dp_j)^2\\
&~\geq \frac{\gjt \Dp_j}{\Dt_j}\cdot \Dp_j - \frac{1}{\alpha}\cdot \frac{\gjt (\Dp_j)^2}{\Dt_j} - |g_j - \tg_j| \cdot |\Dp_j|\comm{as $\Dt_j\leq 1$}\\
&~= \left(1-\frac{1}{\alpha}\right)\frac{\gjt (\Dp_j)^2}{\Dt_j} - |g_j - \tg_j| \cdot |\Dp_j|.
\end{align*}
\end{pfof}

\begin{pfof}{Lemma \ref{lem:error-of-gradient}}
\newcommand{\tpM}{\tp_{\max}}
\newcommand{\tpm}{\tp_{\min}}
\newcommand{\rp}{\mathring{p}}
We begin by showing
\begin{equation}\label{eq:error-of-gradient-1}
\tgjM - \tgjm \leq 2 \sum_{i=1}^m \UH{k_i}{j}{\beta_i}{t}{p_j^t} \cdot |\Dp_{k_i}|.
\end{equation}

First of all, we define a few useful notations.
Let $\tpM$ and $\tpm$, respectively, denote the $\tp$-values at which $\nabla_j \phi(\tp)$ yields $\tgjM$ and $\tgjm$.
Let $p^{(t_1,t]}_{k,\min} := \min_{t'\in (t_1,t]} p^{t'}_k$ and $p^{(t_1,t]}_{k,\max} := \max_{t'\in (t_1,t]} p^{t'}_k$.
Let $\beta_0 := \tau_j$.

To prove \eqref{eq:error-of-gradient-1}, we first construct a path $P$ that connects $\tpM$ and $\tpm$,
with each edge in $P$ corresponding to a price update between times $\tau_j$ and $t$.
The construction builds two paths, $P^s$, starting at $\tpM$, and $P^e$, starting at $\tpm$.
Note that $\tpM,\tpm \in \PP{j}{\tau_j}{t}{p_j^t}$,
and for all $k\neq j$, $(\tpM)_k,(\tpm)_k \in \left[p^{(\beta_0,t]}_{k,\min},p^{(\beta_0,t]}_{k,\max}\right]$.
$P^s$ and $P^e$ will be constructed in $m$ steps that correspond to the $m$ price updates at times $\beta_1,\beta_2,\cdots,\beta_m$.
By the end of the $\ell$-th step, our construction ensures that the end points of $P^s$ and $P^e$ are in the set $\PP{j}{\beta_\ell}{t}{p_j^t}$.
Hence, by the end of the $m$-th step, the end points of $P^s$ and $P^e$ are in the set $\PP{j}{\beta_m}{t}{p_j^t}$, which is a singleton, so the two end points must be equal.
This allows $P^s$ and $P^e$ to be concatenated at their end points to form the path $P$.
The specifics of the construction are as follows:
\begin{enumerate}
\item
Let $\rp^s$ and $\rp^e$, respectively, denote the end points of $P^s$ and $P^e$, i.e.~
initially, $\rp^s = \tpM$ and $\rp^e = \tpm$.
\item
For $i=1\cdots m$, do:
\newcommand{\Li}{p^{(\beta_i,t]}_{k_i,\min}}
\newcommand{\Ri}{p^{(\beta_i,t]}_{k_i,\max}}
\begin{itemize}
\item
Suppose $\Span\left\{\rp^s_{k_i},\rp^e_{k_i}\right\} = [l_i,r_i]$.
WLOG, suppose that $\rp^s_{k_i} = l_i$.\footnote{If $\rp^e_{k_i} = l_i$, swap the roles of $P^s$ and $P^e$ in the current for loop.}

Note that by the end of the last step, the construction ensures that
$l_i,r_i \in \left[p^{(\beta_{i-1},t]}_{k_i,\min},p^{(\beta_{i-1},t]}_{k_i,\max}\right]$.

Also, note that at most one of the strict inequalities $\Li > p^{(\beta_{i-1},t]}_{k_i,\min}$ and
$\Ri < p^{(\beta_{i-1},t]}_{k_i,\max}$ holds,
and hence $l_i < \Li <  \Ri < r_i$ is not possible.

\item For any $p$, let $p' = (p_{-k},x)$ be the vector such that $p'_k = x$, and for all $h \neq k$, $p'_h = p_h$.

Depending on the values of $l_i,r_i,\Li,\Ri$, there are five cases.
\begin{enumerate}
\item If $\Li \leq l_i \leq r_i \leq \Ri$, do nothing.
\item If $l_i < \Li \leq r_i \leq \Ri$, let $\rp ' = \left(\rp^s_{-k_i},\Li\right)$;
in $P^s$, connect $\rp^s$ to $\rp '$, and update $\rp^s$ to $\rp '$.
\item If $l_i \leq r_i < \Li \leq \Ri$,\\
-~let $\rp ' = \left(\rp^s_{-k_i},\Li\right)$; in $P^s$, connect $\rp^s$ to $\rp '$, and update $\rp^s$ to $\rp '$.\\
-~let $\rp '' = \left(\rp^e_{-k_i},\Li\right)$; in $P^e$, connect $\rp^e$ to $\rp ''$, and update $\rp^e$ to $\rp ''$.
\item If $\Li \leq l_i \leq \Ri < r_i$, let $\rp ' = \left(\rp^e_{-k_i},\Ri\right)$;
in $P^e$, connect $\rp^e$ to $\rp '$, and update $\rp^e$ to $\rp '$.
\item If $\Li \leq \Ri < l_i \leq r_i$,\\
-~let $\rp ' = \left(\rp^s_{-k_i},\Ri\right)$; in $P^s$, connect $\rp^s$ to $\rp '$, and update $\rp^s$ to $\rp '$.\\
-~let $\rp '' = \left(\rp^e_{-k_i},\Ri\right)$; in $P^e$, connect $\rp^e$ to $\rp ''$, and update $\rp^e$ to $\rp ''$.
\end{enumerate}
\end{itemize}
\item Concatenate $P^s$ and $P^e$ at $\rp^s = \rp^e$ to form the path $P$.
\end{enumerate}

There are at most $2m$ edges in the path $P$, with at most two edges added in each of the $m$ steps.
Note that the length of each edge added in the $i$-th step is at most $|\Dp_{k_i}|$,
so by simple calculus,
the change to $\nabla_j(p')$ along each such edge
is at most $\UH{k_i}{j}{\beta_i}{t}{p_j^t} \cdot |\Dp_{k_i}|$.
This yields \eqref{eq:error-of-gradient-1}.

\bigskip

To prove \eqref{eq:error-of-gradient-2} and \eqref{eq:error-of-gradient-3},
first note that since $\PP{j}{\beta_i}{t}{p_j^t} \subset \PP{j}{\tau_j}{t}{p_j^t}$,
$\UH{k_i}{j}{\beta_i}{t}{p_j^t} \leq \UH{k_i}{j}{\tau_j}{t}{p_j^t}$.
\begin{align}
\text{Then}~~~~~\sum_{i=1}^m \frac{1}{\eta_i} \UH{k_i}{j}{\beta_i}{t}{p_j^t} \Dt_{k_i}
\leq \sum_{i=1}^m \frac{1}{\eta_i} \UH{k_i}{j}{\tau_j}{t}{p_j^t} \Dt_{k_i}
&\leq \sum_{k\neq j} \frac{1}{\bareta_k} \UH{k}{j}{\tau_j}{t}{p_j^t} \sum_{i:k_i = k} \Dt_{k_i}~~~~~~~~~~~~~\nonumber\\
&\leq 2\sum_{k\neq j} \frac{1}{\bareta_k} \UH{k}{j}{\tau_j}{t}{p_j^t}.\label{eq:lem-error-gradient-inter}
\end{align}
The last inequality holds since $\sum_{i: k_i = k} \Dt_{k_i} \leq 1 + (t-\tau_j) \leq 2$.

\medskip

The proof of \eqref{eq:error-of-gradient-2}:
\begin{align*}
|\mu| \cdot \left(\tgjM - \tgjm\right) &\leq 2 \sum_{i=1}^m \UH{k_i}{j}{\beta_i}{t}{p_j^t} \cdot |\Dp_{k_i}|\cdot |\mu|\comm{by Eqn.~\eqref{eq:error-of-gradient-1}}\\
&\leq \sum_{i=1}^m \UH{k_i}{j}{\beta_i}{t}{p_j^t} \cdot \left[\frac{\mu^2 \Dt_{k_i}}{\eta_i} + \frac{\eta_i (\Dp_{k_i})^2}{\Dt_{k_i}}\right]\comm{AM-GM ineq.}\\
&\leq 2 \mu^2 \sum_{k\neq j} \frac{1}{\bareta_k} \UH{k}{j}{\tau_j}{t}{p_j^t} + \sum_{i=1}^m \eta_i \cdot  \UH{k_i}{j}{\beta_i}{t}{p_j^t}\frac{(\Dp_{k_i})^2}{\Dt_{k_i}}.\comm{by Eqn.~\eqref{eq:lem-error-gradient-inter}}
\end{align*}

The proof of \eqref{eq:error-of-gradient-3}:
\begin{align*}
& \left(\tgjM - \tgjm\right)^2\\
&\leq 4 \sum_{\ione=1}^m \sum_{\itwo=1}^m \UH{k_\ione}{j}{\beta_\ione}{t}{p_j^t}\cdot \UH{k_\itwo}{j}{\beta_\itwo}{t}{p_j^t}\cdot \left|\Dp_{k_\ione}\right|
\cdot \left|\Dp_{k_\itwo}\right|\comm{by Eqn.~\eqref{eq:error-of-gradient-1}}\\
&\leq 2\sum_{\ione=1}^m \sum_{\itwo=1}^m \UH{k_\ione}{j}{\beta_\ione}{t}{p_j^t} \cdot \UH{k_\itwo}{j}{\beta_\itwo}{t}{p_j^t} \cdot \left[\frac{\left(\Dp_{k_\ione}\right)^2\eta_\ione \Dt_{k_\itwo}}{\eta_\itwo \Dt_{k_\ione}} + \frac{\left(\Dp_{k_\itwo}\right)^2\eta_\itwo \Dt_{k_\ione}}{\eta_\ione \Dt_{k_\itwo}}\right]\comm{AM-GM ineq.}\\
&= 2 \sum_{\ione=1}^m \sum_{\itwo=1}^m \UH{k_\ione}{j}{\beta_\ione}{t}{p_j^t} \cdot \UH{k_\itwo}{j}{\beta_\itwo}{t}{p_j^t} \cdot \frac{\left(\Dp_{k_\ione}\right)^2\eta_\ione \Dt_{k_\itwo}}{\eta_\itwo \Dt_{k_\ione}} \\
&\qquad + 2 \sum_{\itwo=1}^m \sum_{\ione=1}^m \UH{k_\itwo}{j}{\beta_\itwo}{t}{p_j^t} \cdot \UH{k_\ione}{j}{\beta_\ione}{t}{p_j^t} \cdot \frac{\left(\Dp_{k_\ione}\right)^2\eta_\ione \Dt_{k_\itwo}}{\eta_\itwo \Dt_{k_\ione}}\\
&\hspace*{1.75in}\comm{swap the indices $i_1$ and $i_2$ in the second double-summation}\\
&= 4 \sum_{\ione=1}^m \sum_{\itwo=1}^m \UH{k_\ione}{j}{\beta_\ione}{t}{p_j^t} \cdot \UH{k_\itwo}{j}{\beta_\itwo}{t}{p_j^t} \cdot
\frac{\left(\Dp_{k_\ione}\right)^2\eta_\ione \Dt_{k_\itwo}}{\eta_\itwo \Dt_{k_\ione}}\\
&= 4 \left(\sum_{\ione=1}^m \eta_\ione \cdot \UH{k_\ione}{j}{\beta_\ione}{t}{p_j^t} \cdot \frac{\left(\Dp_{k_\ione}\right)^2}{\Dt_{k_\ione}}\right)
\left(\sum_{\itwo=1}^m \frac{1}{\eta_\itwo} \UH{k_\itwo}{j}{\beta_\itwo}{t}{p_j^t} \Dt_{k_\itwo}\right)\\
&\leq 8 \left(\sum_{i=1}^m \eta_i \cdot \UH{k_i}{j}{\beta_i}{t}{p_j^t} \cdot \frac{\left(\Dp_{k_i}\right)^2}{\Dt_{k_i}}\right)\left(\sum_{k\neq j} \frac{1}{\bareta_k} \UH{k}{j}{\tau_j}{t}{p_j^t}\right).\comm{by Eqn.~\eqref{eq:lem-error-gradient-inter}}
\end{align*}
\end{pfof}

\begin{pfof}{Lemma \ref{lem:relating-Phi-phi}}
First, we bound the integral terms in $\Phi(p^t,t,\tau)$ (see Eqn.~\eqref{eq:PF}).
Following the derivations of \eqref{eq:bound-integral-1} and \eqref{eq:bound-integral-2}, with $\tg_j$ replaced by $g_j$, yields
$$c_1 \int_{\tau_j}^t \frac{(g_j(t'))^2}{\bgam_j}\,dt' \leq c_1(1+4\eB)\frac{(g_j)^2 \Dt_j}{\bgam_j} + c_1(2+8\eB) \sum_{i=1}^m \xi_j^{\beta_i} \cdot \UH{k_i}{j}{\beta_i}{t}{p_j^t} \frac{(\Dp_{k_i})^2}{\Dt_{k_i}}$$
and hence
$$\sum_j c_1 \int_{\tau_j}^t \frac{(g_j(t'))^2}{\bgam_j}\,dt' \leq c_1(1+4\eB) \sum_j \frac{(g_j)^2 \Dt_j}{\bgam_j} + c_1(2+8\eB) \sum_j \sum_i \xi_j^{\beta_i} \cdot \UH{k_i}{j}{\beta_i}{t}{p_j^t} \frac{(\Dp_{k_i})^2}{\Dt_{k_i}}.$$
When $2-c_2 \geq c_1(2+8\eB)$, as $p_j^t = p_j^{\tau_j+}$,
the double summation in the above inequality is no larger than the double summation in $\Phi(p^t,t,\tau)$.
Thus $\Phi(p^t,t,\tau) \geq \phi(p^t) - c_1(1+4\eB) \sum_j \frac{(g_j)^2 \Dt_j}{\bgam_j}$.

Next, we bound the sum $\sum_j \frac{(g_j)^2 \Dt_j}{\bgam_j}$.
Suppose there are hypothetical updates to all the coordinates at time $t$,
and $p_j$ is updated with the most up-to-date gradient $\tg_j = g_j$ and step size $1/\gamma_j$.
By Lemma \ref{lem:progress-of-phi} and Condition A2,
$\phi^- - \phi^+ \geq \frac{1}{2}\sum_j \frac{(g_j)^2 \Dt_j}{\gamma_j} \ge \frac{1}{2}\sum_j \frac{(g_j)^2 \Dt_j} {\bgam_j} $.
Here $\phi^- = \phi(p^t)$. Thus
$\phi^- - \phi^+ \leq \phi(p^t) - \phi^* = \phi(p^t)$, and hence $\sum_j \frac{(g_j)^2 \Dt_j}{\bgam_j} \leq 2\phi(p^t)$.
\end{pfof}

\newpage

\section{Leontief Fisher Markets}

\begin{lemma}\label{lem:one-big-move-implies-big-progress}
Let $\tau_j,t$ be the times at which two consecutive updates to $p_j$ occur.
If $\gjt$ is controlled and $c_2 \leq 1$, then
$\Phi^{\tau_j+} - \Phi^{t+} \geq \left(1-\frac{1}{\alpha} - 2\eB - 2\eF\right)\frac{\gjt (\Dp_j)^2}{\Dt_j}$.
\end{lemma}
\begin{pf}
This lemma can be proved by slightly modifying the proof of Lemma \ref{lem:Phi-at-update};
we will use the notations defined therein.

By Lemma \ref{lem:Phi-at-update}, $\Phi$ does not increase at the updates made in the time interval $(\tau_j,t)$.
By \eqref{eq:cont-progress},
$$\Phi^{\tau_j+} - \Phi^{t-} \geq c_1 \int_{\tau_j}^t \frac{(g_j(t'))^2}{\bgam_j}\,dt' = E_2.$$
By \eqref{eq:change_of_Phi_inter},
$$\Phi^{t-} - \Phi^{t+} \geq \left(1-\frac{1}{\alpha}\right)\frac{\gjt (\Dp_j)^2}{\Dt_j} - E_1 - E_2 + (2-c_2)\sum_i \xi^{\beta_i}_j \cdot \UH{k_i}{j}{\beta_i}{t}{p_j^{\tau_j+}} \frac{(\Dp_{k_i})^2}{\Dt_{k_i}} - E_3.$$

Combining the two inequalities above yields
\begin{align*}
& \Phi^{\tau_j+} - \Phi^{t+} = \left(\Phi^{\tau_j+} - \Phi^{t-}\right) + \left(\Phi^{t-} - \Phi^{t+}\right)\\
&\geq \left(1-\frac{1}{\alpha}\right)\frac{\gjt (\Dp_j)^2}{\Dt_j} + (2-c_2)\sum_i \xi^{\beta_i}_j \cdot \UH{k_i}{j}{\beta_i}{t}{p_j^{\tau_j+}} \frac{(\Dp_{k_i})^2}{\Dt_{k_i}}  - E_1 - E_3.
\end{align*}
The result follows on noting that $p_j^{\tau_j+} = p_j^t$ and
by applying the bounds on $E_1$ and $E_3$ in the proof of Lemma \ref{lem:Phi-at-update}.
\end{pf}

Let $U=\max\left\{\max_j\{\pc_j\}, 2 \sum_i e_i\right\}$ be an upper bound on the prices throughout the tatonnement process \cite{CCD2013}.

\begin{lemma}\label{lem:big-moves-in-one-time-unit-implies-big-overall-progress}
Let $\delta = 1-\frac{1}{\alpha} - 2\eB - 2\eF$.
Suppose that there are consecutive updates to $p_j$ at times $\Upsilon_0 < \Upsilon_1 < \cdots < \Upsilon_m$, where $\Upsilon_m - \Upsilon_0 \leq 2$.
If $\left|p_j^{\Upsilon_0+} - p_j^{\Upsilon_m+}\right| \geq \epsilon$, where $\epsilon \leq 1$, then $\Phi^{\Upsilon_0+} - \Phi^{\Upsilon_m+} \geq \delta \epsilon^2\cdot \min\left\{\frac{1}{2},\frac{1}{3\lambda U}\right\}$.
\end{lemma}
\begin{pf}
\newcommand{\tz}{\tilde{z}}
For $q=1,2,\cdots,m$, let $\Delta\pjq$ be the change to $p_j$ at the update timed $\Upsilon_q$, and let $\tzjq$ be the $\tz$-value used for the update,
i.e.~$\gamma_j^{\Upsilon_q} = \frac{\max\{1,\tzjq\}}{\lambda p_j^{\Upsilon_q}}$ and
$\Delta\pjq = \lambda p_j^{\Upsilon_q}\cdot \min\{1,\tzjq\}\cdot \Dt_q$.

We will use Lemma \ref{lem:one-big-move-implies-big-progress} to give a lower bound on the decrease to $\Phi$ between times $\Upsilon_0+$ and $\Upsilon_m+$.
If $\tzjq < 1$, then
$$\frac{\gamma_j^{\Upsilon_q}(\D\pjq)^2}{\Dt_q} = \frac{1}{\lambda p_j^{\Upsilon_q}}\frac{(\D\pjq)^2}{\Dt_q}\geq \frac{1}{\lambda U}\frac{(\D\pjq)^2}{\Dt_q}.$$
If $\tzjq \geq 1$, then
$$\frac{\gamma_j^{\Upsilon_q}(\D\pjq)^2}{\Dt_q} = \frac{\tzjq}{\lambda p_j^{\Upsilon_q}}\cdot \lambda^2 \left(p_j^{\Upsilon_q}\right)^2 \Dt_q = \lambda p_j^{\Upsilon_q} \tzjq \Dt_q \geq |\D\pjq|.$$
By Lemma \ref{lem:one-big-move-implies-big-progress},
\begin{align*}
\Phi^{\Upsilon_0+} - \Phi^{\Upsilon_m+} = \sum_{q=1}^m \left(\Phi^{\Upsilon_{q-1}+} - \Phi^{\Upsilon_q+}\right) &\geq \delta \sum_{q=1}^m \frac{\gamma_j^{\Upsilon_q} (\D\pjq)^2}{\Dt_q}\\
&\geq \frac{\delta}{\lambda U}\sum_{q:\tzjq < 1} \frac{(\D\pjq)^2}{\Dt_q} + \delta \sum_{q:\tzjq\geq 1} |\D\pjq|.
\end{align*}

By the assumption $|p_j^{\Upsilon_0+} - p_j^{\Upsilon_m+}| \geq \epsilon$, $\sum_{q=1}^m |\D\pjq| \geq \epsilon$.
Let $\sigma := \epsilon^{-1} \sum_{q:\tzjq\geq 1} |\D\pjq|$. Then $\sum_{q:\tzjq<1} |\D\pjq| \geq \max\{0,(1-\sigma)\epsilon\}$.
By the Cauchy-Schwarz inequality,
\begin{align*}
\left[\max\{0,(1-\sigma)\epsilon\}\right]^2 \leq \left(\sum_{q:\tzjq<1} |\D\pjq|\right)^2 &= \left(\sum_{q:\tzjq<1} \left|\frac{\D\pjq}{\sqrt{\Dt_q}}\right|\cdot \sqrt{\Dt_q}\right)^2\\
&\leq \left(\sum_{q:\tzjq<1} \frac{(\D\pjq)^2}{\Dt_q}\right)\left(\sum_{q:\tzjq<1} \Dt_q\right)\\
&\leq 3\sum_{q:\tzjq<1} \frac{(\D\pjq)^2}{\Dt_q},
\end{align*}
i.e.~$\sum_{q:\tzjq<1} \frac{(\D\pjq)^2}{\Dt_q} \geq \frac{1}{3}\left[\max\{0,(1-\sigma)\epsilon\}\right]^2$. Then
$$\Phi^{\Upsilon_0+} - \Phi^{\Upsilon_m+} \geq \frac{\delta}{3\lambda U} \left[\max\{0,(1-\sigma)\epsilon\}\right]^2 + \delta \sigma \epsilon.$$
The minimum value of the right hand side is at least $\delta \epsilon^2\cdot \min\left\{\frac{1}{2},\frac{1}{3\lambda U}\right\}$.
\end{pf}

\begin{cor}\label{cor:small-move-in-one-time-unit}
For any $\epsilon > 0$, there exists a finite time $T_\epsilon$ such that for any good $j$,
any $t\geq T_\epsilon$, and any $0\leq \Dt \leq 1$, $|p_j^t - p_j^{t+\Dt}| \leq \epsilon$.
\end{cor}

\begin{pfof}{Theorem \ref{thm:CES-Fisher-cnvge} for the Leontief case}
The proof comprises four steps. We need the following definitions:
for any two price vectors $p^A$ and $p^B$, let $d(p^A,p^B)$ denote the $L_1$ norm distance between the two price vectors,
i.e.~$d(p^A,p^B) = \sum_j |p^A_j - p^B_j|$.
For any two sets of price vectors $P^A$ and $P^B$, let $d(P^A,P^B) := \inf_{p^A\in P^A,p^B\in P^B} d(p^A,p^B)$.

\medskip

\noindent\textbf{Step 1.} Let $\Omega$ be the set of limit points of a tatonnement process. We show that $\Omega$ is non-empty and connected.

\smallskip

Since all prices remain bounded by $U$ throughout the tatonnement process, $\Omega$ is non-empty.

Suppose $\Omega$ is not connected. Let $\Omega_a$ denote a connected component of $\Omega$, and let $\Omega_b = \Omega\setminus \Omega_a$.
Suppose $d(\Omega_a,\Omega_b) = \ep ' > 0$.
By the definition of limit points, there exists a finite time $\Upsilon_{\epsilon'}$ such that thereafter the prices in the tatonnement process
are always within an $\epsilon'/4$-neighborhood of either $\Omega_a$ or $\Omega_b$.
This forces an infinite number of updates, each separated by at least one time unit, such that each update makes a change to a price by at least at least $\epsilon'/2$.
This contradicts Corollary \ref{cor:small-move-in-one-time-unit}.

\medskip

\noindent\textbf{Step 2.} Recall that a market equilibrium is a price vector $\ps$ at which for each $j$, $\ps_j>0$ implies $z_j(\ps)=0$ and $\ps_j=0$ implies $z_j(\ps)\leq 0$. We define a pseudo-equilibrium: a price vector $\tp$ is a pseudo-equilibrium if for each $j$, $\tp_j>0$ implies $z_j(\tp)=0$. Note that every market equilibrium is a pseudo-equilibrium. We show that all limit points in $\Omega$ are pseudo-equilibria.

\smallskip

\newcommand{\ringp}{\grave{p}}
Suppose not. Let $p'\in \Omega$ be a price vector which is not a pseudo-equilibrium, i.e.~there exists $j$ such that $p'_j>0$ but $z_j(p')\neq 0$.
Let $\epsilon$ be a positive number such that for any price vector $\ringp$ in the $\ep$-neighborhood of $p'$,
$\ringp_j\geq p'_j/2$ and $|z_j(\ringp)| \geq |z_j(p')/2|$.
By the definition of limit points, the tatonnement process enters the $(\ep/2)$-neighborhood of $p'$ infinitely often.
By Corollary \ref{cor:small-move-in-one-time-unit}, there exists a finite time such that subsequently,
every time the tatonnement process enters the $\epsilon/2$-neighborhood of $p'$, it stays in the $\epsilon$-neighborhood of $p'$ for at least one time unit.
By Eqn.~\eqref{eq:cont-progress}, $\Phi$ drops by at least $\lambda (p'_j/2) (z_j(p')/2)^2$ during each such stay in the $\epsilon$-neighborhood of $p'$.
This is a contradiction since $\Phi$ is positive throughout and hence cannot drop by at least $\lambda (p'_j/2) (z_j(p')/2)^2$ infinitely often.

\medskip

\noindent\textbf{Step 3.} We show that the excess demands at all limit points in $\Omega$ are identical.

\smallskip

For every subset of goods $S$, let $\Omega_S = \{p'\in \Omega\,|\,p'_k > 0 \Leftrightarrow k\in S\}$.
For each buyer, there are two cases:
\begin{itemize}
\item \textbf{if the buyer wants at least one good in $S$, say good $\ell$:}\\
Observe that by the definition of pseudo-equilibrium and Step 2, every price vector in $\Omega_S$,
excluding the zero prices in the price vector,
is a market equilibrium for the sub-Leontief-market comprising the goods in $S$.
Codenotti and Varadarajan \cite{CV2004} pointed out that the demands for the goods in $S$ of each buyer are identical at every market equilibrium of the sub-Leontief market,
and hence also in the original Leontief market.
So the buyer demands the same positive but finite amount of good $\ell$ at every price vector in $\Omega_S$
in the original market.
Also note that the buyer always demands the goods in the original market in a fixed proportion.
This forces the demands for the goods not in $S$ of the buyer are also identical at every price vector in $\Omega_S$.
\item \textbf{if the buyer wants no good in $S$:}\\
Then the buyer demands infinite amount of each good that she wants, and demands zero amount of each good that she does not want.
\end{itemize}
In either case, the buyer's demands for each good at every price vector in $\Omega_S$ are identical,
and hence also the total demand for each good.

Then consider a graph $G$ with each vertex corresponding to a subset of goods $S$ such that $\Omega_S$ is non-empty,
and two vertices $S_1,S_2$ being adjacent if and only if $d\left(\Omega_{S_1},\Omega_{S_2}\right) = 0$.
Since excess demands are a continuous function\footnote{The range of the excess demand functions is the extended real line $\rr \cup \{+\infty\}$; continuity of the excess demand function is w.r.t.~the usual topology on the extended real line. To be specific, if $z_k(p)=+\infty$ for some $p$ and $k$, then for any $M\in \rr$, there exists an $\ep_M>0$ such that $z_k(p)\geq M$ in the $\ep_M$-neighborhood of $p$.} of prices,
if $S_1$ and $S_2$ are adjacent, then the excess demands for all goods at every price vector in $S_1\cup S_2$ are identical.
By Step 1, the graph $G$ is connected, thus the excess demands at all limit points in $\Omega$ are identical.

\bigskip

\noindent\textbf{Step 4.} We show that every limit point in $\Omega$ is indeed a market equilibrium.

\smallskip

Suppose not, i.e.~there exists a limit point $p'$ in $\Omega$ which is a pseudo-equilibrium but not a market equilibrium,
i.e.~there exists $k$ such that $p'_k=0$ but $z_k(p')>0$.
By Step 3, $z_k$ is positive at every limit point in $\Omega$, and hence every $p_k$ at every limit point must be zero.
By the definition of limit points, for any $\epsilon>0$, beyond a finite time, the tatonnement process must stay within the $\epsilon$-neighborhood of $\Omega$ thereafter.
By choosing a sufficiently small $\epsilon$, $z_k$ is bounded away from zero in the $\epsilon$-neighborhood of $\Omega$, and hence $p_k$ increase indefinitely and eventually $p_k$ becomes so large that the tatonnement process must leave the $\epsilon$-neighborhood of $\Omega$, a contradiction.
\end{pfof}

\newpage

\section{Ongoing Complementary-CES Fisher Markets}

\newcommand{\tz}{\tilde{z}}
\newcommand{\kj}{\kappa_j}
\newcommand{\kv}{\kj v_j}
\newcommand{\lj}{\lambda_j}
\newcommand{\bz}{\bar{z}}
\newcommand{\wh}{\mathcal{W}}

The tatonnement process which we described in Section \ref{sect:CES} is a two-stage process.
In the first stage, the buyers repeatedly \emph{report} their demands to sellers according to the current prices,
then the sellers update the prices with the reported demands.
The first stage continues until the market reaches a market equilibrium, 
and then trades occur in the second stage.
Clearly, this is not a plausible real-world market dynamic.

In order to have a more realistic setting for a price adjustment algorithm,
it would appear that out-of-equilibrium trade must be allowed,
so as to generate the demand imbalances that then induce price adjustments.
In an attempt to build a more realistic market model,
Cole and Fleischer \cite{CF2008} introduced the \emph{Ongoing market model}.
In an ongoing Fisher market, the market repeats over an unbounded number of time intervals called \emph{days}.
Each day, the seller of each good receives one new unit of the good, and each buyer $i$ is given $e_i$ amount of money.
In that day, each buyer $i$ purchases a utility-maximizing bundle of goods of cost at most $e_i$.

But then there needs to be a way for seller to handle excess supply/demand.
To this end, for each good $j$ there is a \emph{warehouse} of finite capacity $\chi_j$ which can meet excess demand and store excess supply.
When there is surplus (supply exceeds demand),
it is stored in the warehouse;
when there is excess demand (demand exceeds supply),
good is taken from the warehouse to meet the excess demand.
The sellers change prices as needed to ensure their warehouses neither overfill nor run out of goods.

Given initial prices $p^0$, initial warehouses stocks $v^0$, where $0<v^0_j<\chi_j$ for each good $j$,
and \emph{ideal warehouse stocks} $v^*$,
the task is to repeatedly adjust prices so as to converge to a market equilibrium
with the warehouse stocks converging to their ideal values;
for simplicity, we suppose that $v^*_j = \chi_j/2$ for each good $j$.
$v_j$ will denote the difference between the content of the warehouse of good $j$ and $v^*_j$; hence $v_j \in [-\chi_j/2,\chi_j/2]$.

In an ongoing Fisher market, the sellers adjust the prices of their goods.
In order to have progress, the sellers are required to update prices at least once per day.
However, there is no upper bound on the frequency of price changes.
This entails measuring demand on a finer scale than day units.
Accordingly, we assume that each buyer spends their money at a uniform rate throughout the day,
and hence \emph{instantaneous demand} and \emph{instantaneous excess demand} for good $j$
at any time $t\in \rr^+$ can be readily defined;
we denote them by $x_j^t$ and $z_j^t$ respectively.

In this section, we analyse ongoing complementary-CES Fisher markets.
Recall that for a complementary-CES Fisher market, tatonnement is equivalent to gradient descent
on the convex function
$\phi(p) = \sum_j p_j + \sum_i \hat{u}_i(p)$,
where $\hat{u}_i(p)$ is the optimal utility that buyer $i$ can attain at prices $p$.
We will introduce new potential functions, which incorporate $\phi$ as a component,
for the ongoing market analysis.

We use the following price update rule, which is a variant of \eqref{eq:async-tat-CES-rule},
and which ensures convergence to the ideal warehouse stocks as well as to the market equilibrium:
\begin{equation}\label{eq:update-rule-with-warehouse}
p_j' = p_j \cdot \left(1 + \lj \cdot \min\{\tz_j - \kv,1\} \cdot \Dt_j\right),
\end{equation}
where $\lj,\kj$ are small constants. Note that $\gjt = \frac{1}{p_j \lj}\cdot \max\{1,\tz_j-\kv\}$.

\begin{theorem}\label{thm:tat-with-warehouse}
If $\lj \leq 1/60$ for all $j$, then there exists $\kj>0$ such that price updates using Rule \eqref{eq:update-rule-with-warehouse} converge toward the market equilibrium in any complementary-CES Fisher market,
with the warehouse stocks converging to their ideal values.
\end{theorem}

First, we impose the following bounds on $\lj$ and $\kj$.
\begin{enumerate}
\item[B1.] $\lj \leq 1/60$;
\item[B2.] $\kj/\lj \leq 1/10$ (this, together with Condition B1, yields $\kj \leq 1/600$);
\item[B3.] $|\kv| \leq 1/10$ always (such $\kj$ exist since the warehouse sizes are bounded).
\end{enumerate}

We will impose more bounds on $\kj$,
but eventually we will show that, given any fixed $\lj$ satisfying Condition B1,
for all $j$, there exist positive $\kj$ that satisfy all these bounds.

We need to be cautious with Condition B3, and also Condtion B4 which we will state later.
At this point, it is not clear that $v_j$ remains bounded throughout the tatonnement process,
so the two conditions might cease to hold no matter how small $\kj$ is set.
We show that this never happens in Section \ref{sect:bounding-warehouse}.

\noindent\textbf{Notations} Let $f\geq 1$. A price vector $p$ is $f$-bounded if, for all $j$,
$\frac{1}{f} \leq \frac{p_j}{\ps_j} \leq f$.
Let $R(f)$ denote the set of all $f$-bounded price vectors.

Our analysis comprises two phases.
Phase 1 finishes when prices are guaranteed to be $1.9$-bounded thereafter,
and then we proceed to Phase 2.
We outline the analysis of the two phases in Sections~\ref{sect:phase1} and~\ref{sect:phase2}, respectively.
We defer most proofs to Section \ref{sect:ongoing-missing-proof}.

One component of the potential functions we will use is (similar to) $\Phi$ as defined in \eqref{eq:PF},
and we will use some results from Sections \ref{sect:analysis} and \ref{sect:CES}.
We deduce the values of $\eB,\eF$ that satisfy Conditions A3 and A4.
Recall that by Property 3 of complementary-CES markets (see the appendix on tatonnement),
if $\frac{\Dp_j}{p_j} \leq 1/6$,
then $\phi(p+\Dp) - \phi(p) - \nabla_j\phi(p) \cdot \Dp_j \leq \frac{1.5 x_j}{p_j} (\Dp_j)^2$,
where $x_j = z_j + 1$.
Let $\tilde{x}_j = \tz_j + 1$.
Recall that $x_j \leq \frac{100}{81} \tilde{x}_j \leq 1.24 \tilde{x}_j$.
Then
\begin{equation}\label{eq:progress-of-phi-inter}
\frac{1.5 x_j}{p_j}\cdot \frac{1}{\gjt} \leq \frac{1.86 \tilde{x}_j}{p_j} \cdot \frac{p_j \lj}{\max\{\tz_j-\kv,1\}}\\
\leq 1.86 \lj \cdot \frac{\tz_j + 1}{\max\{\tz_j-0.1,1\}} \leq \frac{1.86}{60} \cdot 2.1 < \frac{1}{15} < \frac{1}{2},
\end{equation}
and hence $\frac{\gjt}{15} > \frac{1.5 x_j}{p_j}$.
By Lemma \ref{lem:ces-bounds} plus Conditions (A3) and (A4),
we can set 
\begin{eqnarray}
\label{eq:eF-bound}
\eF &=& \frac{1.53}{1.5}\cdot \frac{1}{15} ~=~ 0.068\\
\label{eq:eB-bound}
\text{and}~~~ \eB &=& \frac{1.89}{1.5}\cdot \frac{1}{15} ~=~ 0.084.
\end{eqnarray}

\begin{lemma}\label{lem:progress-of-phi-with-warehouse}
Suppose there is an update to $p_j$ at time $t$ according to rule \eqref{eq:update-rule-with-warehouse}.
Suppose that Conditions B1 and B3 hold.
Let $\phi^-$ and $\phi^+$ denote, respectively, the convex function values just before and just after the update.
Let $z_j = -\nabla_j \phi(p^t)$ and $\tz_j \equiv \tz_j(t)$.
Let $\Dp_j$ be the change to $p_j$ made by the update, i.e.~$\Dp_j := \lj p_j \cdot \min\{\tz_j - \kv,1\} \cdot \Dt_j$.
Then
\begin{equation}\label{eq:progress-of-phi-1}
\phi^- - \phi^+ \geq \frac{1}{2} \frac{(\tz_j)^2 \Dt_j}{\gjt} - \frac{1}{2} \frac{(\kv)^2 \Dt_j}{\gjt} - |z_j - \tz_j|\cdot |\Dp_j|
\end{equation}
and
\begin{equation}\label{eq:progress-of-phi-2}
\phi^- - \phi^+ \geq \frac{41}{60} \frac{\gjt (\Dp_j)^2}{\Dt_j} - \frac{(\kv)^2 \Dt_j}{\gjt} - |z_j - \tz_j|\cdot |\Dp_j|.
\end{equation}
\end{lemma}

\newpage

\subsection{Phase 1}\label{sect:phase1}

\newcommand{\gjsj}{\gamma_j^{\sigj}}

For Phase 1, we use the potential function $\Xi_1 \equiv \Xi_1(p^t,v^t,t,\tau)$:
\begin{align}
\Xi_1 &= \phi(p^t) - c_1 \sum_j \int_{\tau_j}^{t} \frac{(z_j(t'))^2}{\gjsj}\,dt' + \sum_j \sum_i \xi^{\beta_i}_j \UH{k_i}{j}{\beta_i}{\sigj}{p_j^{\tau_j+}}\frac{(\Dp_{k_i})^2}{\Dt_{k_i}} [2-c_2(t - \beta_i)]\nonumber\\
&\qquad\qquad + \sum_j \frac{(\kj v_j^t)^2 (t-\tau_j)}{\gjsj}.\label{eq:Xi1}
\end{align}

When there is no update, we show that
\begin{equation}\label{eq:Xi1-cont}
\frac{d \Xi_1}{dt}
\leq -\sum_j (c_1-\kj) \frac{(z_j^t)^2}{\gjsj} + \sum_j (1+\kj) \frac{(\kj v_j^t)^2}{\gjsj} - c_2 \sum_j \sum_i
\xi^{\beta_i}_j \UH{k_i}{j}{\beta_i}{\sigj}{p_j^{\tau_j+}}\frac{(\Dp_{k_i})^2}{\Dt_{k_i}}.
\end{equation}

When there is an update, we show that
\begin{lemma}\label{lem:Xi1-at-update}
Suppose that there is an update to $p_j$ at time $t$. Suppose that Conditions B1 and B3 hold.
Let $\Xi_1^-$ and $\Xi_1^+$, respectively, denote the values of $\Xi_1$ just before and just after the update.
Then
$$\Xi_1^- - \Xi_1^+ \geq \left(\frac{1}{4} - 1.4 c_1\right)\frac{(\tz_j)^2 \Dt_j}{\gjt} + \left(1- c_2 - 2.7 c_1\right)
\sum_{i=1}^m \xi^{\beta_i}_j \UH{k_i}{j}{\beta_i}{t}{p_j^{\tau_j+}} \frac{(\Dp_{k_i})^2}{\Dt_{k_i}}.$$
\end{lemma}
Thus, by setting $c_1=5/28$ and $c_2=1/2$, $\Xi_1$ does not increase at any update.

Since $\phi$ is strongly convex, in the proof of Theorem \ref{thm:main}(b), we show that
$\sum_j \frac{(z_j)^2}{\bgam_j} \geq D_1 \cdot \phi(p^t)$ for some positive constant $D_1 \leq 1/10$.
Let $\psi := \frac{1}{1.001} \cdot \inf_{p'\notin R(1.9)} \phi(p')$.
We impose an additional condition on $\kj$:
\begin{enumerate}
\item[B4.] $\kj$ are sufficiently small such that $\sum_j \frac{(\kv)^2}{\gjt} \leq \frac{1}{26/D_1 + 4} \psi$ always.
\end{enumerate}

\begin{lemma}\label{lem:Xi1-decrease}
If Condition B4 holds and $\Xi_1 \geq \psi/2$, then $\frac{d \Xi_1}{dt} \leq -\Theta(1) \cdot \Xi_1(t)$.
\end{lemma}

\begin{pf}
Let $H(t)$ denote the sum
$\sum_j \sum_{i=1}^m \xi^{\beta_i}_j \UH{k_i}{j}{\beta_i}{t}{p_j^{\tau_j+}} \frac{(\Dp_{k_i})^2}{\Dt_{k_i}}$
at time $t$.
By \eqref{eq:Xi1} and Condition B4,
$$\phi(p^t) + 2H(t) + \frac{1}{26/D_1 + 4} \psi ~\geq ~ \Xi_1(t) ~\geq ~ \psi/2.$$
Hence
\begin{equation}
\label{eqn:H-bound-one}
\phi(p^t) + 2H(t) ~\geq ~ \left(\frac{1}{2} - \frac{1}{26/D_1 + 4}\right) \psi
\end{equation}
\begin{equation}
\label{eqn:H-bound-two}
\mbox{and}~~~~\phi(p^t) + 2H(t) ~\geq ~ \left(1 - \frac{1}{13/D_1 + 2}\right) \Xi_1(t).
\end{equation}

With our choices of $c_1,c_2$ and Condition B4, \eqref{eq:Xi1-cont} yields
\begin{align}
\frac{d \Xi_1}{dt} &\leq -\sum_j \left(\frac{5}{28}-\kj\right) \frac{(z_j^t)^2}{\gjsj} + \sum_j (1+\kj) \frac{(\kj v_j^t)^2}{\gjsj} - \frac{1}{2}  H(t)\nonumber\\
&\leq -\frac{1}{6} \sum_j \frac{(z_j^t)^2}{\gjsj} + \frac{601}{600}\cdot \frac{1}{26/D_1 + 4} \psi  - \frac{1}{2} H(t)\nonumber\\
&\leq -\frac{D_1}{6} \cdot \phi(p^t)  - \frac{1}{2} H(t) + \frac{601}{600}\cdot \frac{1}{26/D_1 + 4} \psi\nonumber\\
&\leq -\frac{D_1}{6}\left(\phi(p^t) + 2H(t)\right) + \frac{\frac{601}{600}\cdot \frac{1}{26/D_1 + 4}}{\frac{1}{2} - \frac{1}{26/D_1 + 4}}\left(\phi(p^t) + 2H(t)\right)
~~~~~~~~~~~\mbox{(by Eqn.~\eqref{eqn:H-bound-one})}\nonumber\\
&\leq -\frac{D_1}{12}\left(\phi(p^t) + 2H(t)\right)\nonumber\\
&\leq -\frac{D_1}{12}\cdot \left(1 - \frac{1}{13/D_1 + 2}\right) \Xi_1(t)
\hspace*{2in}~ \mbox{(by Eqn.~\eqref{eqn:H-bound-two})}\nonumber\\
&\leq -\frac{D_1}{13}\cdot \Xi_1(t).\label{eq:Xi1-cont-concrete}
\end{align}
\end{pf}

\begin{lemma}\label{lem:Xi1-thereafter}
If $\Xi_1(t_1) < \psi/2$ at some time $t_1$, then $\Xi_1(t) \leq \psi/2$ thereafter.
\end{lemma}
\begin{pf}
Suppose the contrary, i.e.~at some time $t_2>t_1$, $\Xi_1(t_2) > \psi/2$.
Let $T_2$ be the collection of all such $t_2$, and let $t'$ be the infimum of $T_2$.
By Lemma \ref{lem:Xi1-at-update} and our choices of $c_1$ and $c_2$,
$\Xi_1$ never increases at an update.
Hence, for $\Xi_1$ to exceed $\psi/2$ after time $t_1$, it must be due to continuous incrementing.
This forces $\Xi_1(t')=\psi/2$ and $\left. \frac{d\Xi_1}{dt}\right|_{t=t'} \geq 0$.
But these contradict Lemma \ref{lem:Xi1-decrease}.
\end{pf}

Following the proof of Lemma~\ref{lem:relating-Phi-phi}, we obtain that
$\Xi_1 \ge \phi(p^t) -2c_1(1 + 8\eB) \phi(p^t)$, and as $c_1(1 + 8\eB) \le \frac 14$,
$\Xi_1 \ge \frac 12 \phi(p^t)$.
Thus if $\Xi_1 \leq \psi/2$,
then $\phi(p^t)/2 \leq \Xi_1 \leq \psi/2$.
This implies $\phi(p^t) < \min_{p'\notin R(1.9)} \phi(p')$ and thus $p^t\in R(1.9)$.
Lemma \ref{lem:Xi1-decrease} shows that $\Xi_1$ decreases linearly
until it drops below $\psi/2$ at some time $t_1$,
and Lemma \ref{lem:Xi1-thereafter} shows that $\Xi_1$ remains below $\psi/2$ thereafter.
Hence, $\forall t\geq t_1$, $p^t\in R(1.9)$ and we proceed to the analysis of Phase 2.

\newpage

\subsection{Phase 2}\label{sect:phase2}

Phase 2 starts when all prices are guaranteed to be $1.9$-bounded thereafter.
Then each demand is between $\frac{1}{1.9}$ and $1.9$ and hence $-0.5 \leq z_j,\tz_j \leq 0.9$.
Since $|\kv|\leq 0.1$ always, in Phase 2 the update rule \eqref{eq:update-rule-with-warehouse} is equivalent to
\begin{equation}\label{eq:update-rule-with-warehouse-p2}
p_j' = p_j \cdot \left(1 + \lj \cdot (\tz_j - \kv) \cdot \Dt_j\right),
\end{equation}
i.e.~$\gjt = \frac{1}{\lj p_j}$.

In this phase, we will use a new potential function $\Xi_2$, which comprises two main components $\Phi$ and $\wh$.
$\Phi$ reflects how far the current prices are from the market equilibrium,
and $\wh$ accounts for the warehouse imbalances.

\subsubsection{Component $\Phi$}

The first component of $\Xi_2$, $\Phi\equiv \Phi(p^t,t,\tau)$, is
\begin{equation}\label{eq:Xi2-Phi}
\Phi = \phi(p^t) - c_1 \sum_j \int_{\tau_j}^{t} \lj p_j (z_j(t'))^2\,dt' + \sum_j \sum_i \xi^{\beta_i}_{j} \UH{k_i}{j}{\beta_i}{\sigj}{p_j^{\tau_j+}} \frac{(\Dp_{k_i})^2}{\Dt_{k_i}}  \left[6-c_2 (t - \beta_i)\right].
\end{equation}

When there is no update, it is straightforward to show that
\begin{equation}\label{eq:phase2-Phi-cont}
\frac{d \Phi}{dt} = -c_1 \sum_j \lj p_j (z_j^t)^2 - c_2 \sum_j \sum_i \xi^{\beta_i}_j \UH{k_i}{j}{\beta_i}{\sigj}{p_j^{\tau_j+}} \frac{(\Dp_{k_i})^2}{\Dt_{k_i}}.
\end{equation}

When there is an update, we show that
\begin{lemma}\label{lem:Phi-at-update-with-warehouse}
Suppose that there is an update to $p_j$ at time $t$. Suppose that Conditions B1 and B3 hold.
Let $\Phi^-$ and $\Phi^+$, respectively, denote the values of $\Phi$ just before and just after the update.
Then
\begin{align*}
\Phi^- - \Phi^+ &\geq \left(\frac{1}{20} - 1.4 c_1\right) \lambda_j p_j (\tz_j)^2 \Dt_j
+ 0.039 \frac{(\Dp_j)^2}{\lj p_j \Dt_j}
- \frac{19}{20}\lj p_j (\kv)^2 \Dt_j \\
&\qquad + \left(5 - c_2 - 2.7 c_1\right) 
\sum_i \xi^{\beta_i}_j \UH{k_i}{j}{\beta_i}{\sigj}{p_j^{\tau_j+}} \frac{(\Dp_{k_i})^2}{\Dt_{k_i}}.
\end{align*}
\end{lemma}

\subsubsection{Component $\wh$}

Let $f_j := \ln (p_j / \ps_j)$. The second component of $\Xi_2$, $\wh \equiv \wh(p^t,v^t,t,\tau)$, is
$$\wh = \sum_j \frac{\kj}{\lj} \ps_j \left(f_j + \lj v_j\right)^2 - c_3 \sum_j \lj \ps_j (\kv)^2 (t-\tau_j) 
+ 2 \sum_j \kj \lj \ps_j \int_{\tau_j}^t v_j(t') z_j(t')\,dt'.$$

When there is no update, we show that for any $R_1\in \rr^+$,
\begin{equation}\label{eq:phase2-W-cont}
\frac{d\wh}{dt} \leq -c_3 \sum_j (1-\kappa_j) \lj \ps_j (\kv^t)^2 + \sum_j (R_1 + c_3 \lj) \kj \ps_j (z_j^t)^2  + \frac{1}{R_1}\sum_j \kj \ps_j (f_j)^2.
\end{equation}
We will choose an appropriate value of $R_1$ at the end.

\begin{lemma}\label{lem:W-at-update-with-warehouse}
Suppose that there is an update to $p_j$ at time $t$.
Suppose that Conditions B1--B3 hold.
Let $\wh^-$ and $\wh^+$, respectively, denote the values of $\wh$ just before and just after the update. Then for any $R_2\in \rr^+$,
\begin{align*}
\wh^- - \wh^+ &\geq \left(0.858 - \frac{c_3}{1.9}\right) \lj p_j (\kv)^2 \Dt_j - 0.0235 \lj p_j (\tz_j)^2 \Dt_j \\
&\qquad\qquad - 3.809 \sum_i \xi^{\beta_i}_j \UH{k_i}{j}{\beta_i}{\sigj}{p_j^{\tau_j+}} \frac{(\Dp_{k_i})^2}{\Dt_{k_i}}\\
&\qquad\qquad - 0.101 \kj \frac{\ps_j (f_j)^2 \Dt_j}{R_2} - 1.92 R_2 \frac{(\Dp_j)^2}{\lj p_j \Dt_j}.
\end{align*}
\end{lemma}
We will choose an appropriate value of $R_2$ at the end.

\subsubsection{Ultimate Potential Function $\Xi_2$}

The ultimate potential function $\Xi_2 \equiv \Xi_2(p^t,v^t,t,\tau)$ is
$$\Xi_2 := \Phi + 1.2 \wh + 0.1212 \sum_j \frac{\kj \ps_j (f_j)^2}{R_2} (t-\tau_j).$$

From Lemmas~\ref{lem:Phi-at-update-with-warehouse} and~\ref{lem:W-at-update-with-warehouse}, we deduce that
\begin{align}
& (\Xi_2)^- - (\Xi_2)^+ \nonumber\\
&\geq \left(0.039 - 2.304 R_2\right) \frac{(\Dp_j)^2}{\lj p_j \Dt_j}
+\left(0.0218 - 1.4 c_1\right) \lj p_j (\tz_j)^2 \Dt_j
+ \left(0.0796 - \frac{12c_3}{19}\right) \lj p_j (\kv)^2 \Dt_j \nonumber\\
&\qquad + \left(0.4292 - c_2 - 2.7 c_1\right) \sum_i \xi^{\beta_i}_j \UH{k_i}{j}{\beta_i}{\sigj}{p_j^{\tau_j+}} \frac{(\Dp_{k_i})^2}{\Dt_{k_i}}.\label{eq:Xi2-at-update}
\end{align}

From \eqref{eq:phase2-Phi-cont}, \eqref{eq:phase2-W-cont}
and the fact that $p_j \ge p^*_j/1.9$, we deduce that
\begin{align}
\frac{d \Xi_2}{dt} & \leq \sum_j \left[\frac{2.28 \kj}{\lj}(R_1 + c_3\lj) - c_1\right] \lj p_j (z_j^t)^2 -1.2 c_3 \sum_j (1-\kappa_j) \lj \ps_j (\kv^t)^2 \nonumber\\
&\qquad - c_2 \sum_j \sum_i \xi^{\beta_i}_j \UH{k_i}{j}{\beta_i}{\sigj}{p_j^{\tau_j+}} \frac{(\Dp_{k_i})^2}{\Dt_{k_i}}
+ \left(\frac{1.2}{R_1} + \frac{0.1212}{R_2}\right)\sum_j \kj \ps_j (f_j)^2.\label{eq:Xi2-cont}
\end{align}

\medskip

We also show the following upper and lower bounds on $\Xi_2$.
\begin{equation}
\label{eq:c2-cond}
\text{If}~~~ 2-c_2 \geq 2.7 c_1, \hspace*{4.7in}
\end{equation}
\begin{equation}\label{eq:Xi2-lower-bound}
\Xi_2 \geq (1- 2.7 c_1) \phi(p^t)
- 1.2 \sum_j \frac{\kj}{\lj} \ps_j (f_j)^2
- 20 \sum_j \kj \lj p_j (z_j)^2
+ \sum_j \left(\frac{1}{5} - 1.2 c_3 \kj\right) \kj \lj \ps_j (v_j)^2.
\end{equation}
Also,
\begin{align}
\Xi_2 &\leq \phi(p^t) + \sum_j \left(\frac{2.4}{\lj} + \frac{0.1212}{R_2}\right) \kj \ps_j (f_j)^2 + 20 \sum_j \kj \lj p_j (z_j)^2 \nonumber\\
&\qquad + 10 \sum_j \sum_i \xi^{\beta_i}_j \UH{k_i}{j}{\beta_i}{\sigj}{p_j^{\tau_j+}} \frac{(\Dp_{k_i})^2}{\Dt_{k_i}}
+ 3.6 \sum_j \kj \lj \ps_j (v_j)^2.\label{eq:Xi2-upper-bound}
\end{align}

\medskip

In the next lemma, we show that $\sum_j \ps_j (f_j)^2 = O(1) \cdot \sum_j p_j (z_j)^2$,
with the hidden constant in $O(1)$ depending on $\max_i \theta_i$,
where $\theta_i$ is the parameter of the CES utility function of buyer $i$.

\begin{lemma}\label{lem:z-to-f}
Let $R := \left\{p'\,\left|\, \forall j,~\frac{1}{1.9} p_j^* \leq p'_j \leq 1.9 p_j^*\right.\right\}$
and $\bar{\theta} = \max_i \theta_i$.For all $p'\in R$,
$$\sum_j p_j^* (f_j)^2 \leq \overline{M}\sum_j p_j '(z_j)^2,$$
where
$\overline{M} = \left(1-\bar{\theta}\right)^{-1}
\max\left\{26.56~,~6.64\bar{\theta}\left(1+\bar{\theta}-2^{\bar{\theta}}\right)^{-1}\right\}$.
\end{lemma}

\medskip

Finally, we choose parameters $R_1,R_2,c_1,c_2,c_3$ such that $\Xi_2$ never increases at an update,
and if there is no update, then $\frac{d \Xi_2}{dt} \leq -\Theta(1)\cdot \Xi_2$.
Set $R_2 = 39/2304$, $c_1 = \frac{0.0218}{1.4} \approx 0.0156$, $c_3 = \frac{19\times 0.0796}{12}\approx 0.1260$,
$c_2 = 0.3855$ and $R_1 = 1$. By choosing sufficiently small $\kj$, \eqref{eq:Xi2-cont} and Lemma \ref{lem:z-to-f} yield
$$\frac{d \Xi_2}{dt} \leq -\Theta(1)\cdot \sum_j \lj p_j (z_j^t)^2 - \Theta\left(\min_j \kj\right)\cdot \sum_j \kj \lj \ps_j (v_j^t)^2 - \Theta(1)\cdot \sum_j \sum_i \xi^{\beta_i}_j \UH{k_i}{j}{\beta_i}{\sigj}{p_j^{\tau_j+}} \frac{(\Dp_{k_i})^2}{\Dt_{k_i}}.$$
Also, by \eqref{eq:Xi2-upper-bound}, Lemma \ref{lem:z-to-f} and the fact that $\phi(p) \leq \Theta(1)\cdot \sum_j p_j (z_j)^2$ \cite[Lemma 6.3]{CCD2013} yield
$$\Xi_2 \leq \Theta(1) \cdot \sum_j p_j (z_j^t)^2 + \Theta(1)\cdot \sum_j \kj \lj \ps_j (v_j^t)^2 + \Theta(1) \cdot \sum_j \sum_i \xi^{\beta_i}_j \UH{k_i}{j}{\beta_i}{\sigj}{p_j^{\tau_j+}} \frac{(\Dp_{k_i})^2}{\Dt_{k_i}}.$$
Thus $\frac{d \Xi_2}{dt} \leq -\Omega\left(\min_j \kj\right)\cdot \Xi_2$.

Further, \eqref{eq:Xi2-lower-bound} and the fact that $\phi(p) \geq \Theta(1)\cdot \sum_j p_j (z_j)^2$ \cite[Lemma 6.2]{CCD2013} yield
\begin{equation}\label{eq:Xi2-lower-concrete}
\Xi_2 \geq \Theta(1)\cdot \phi(p^t) + \Theta(1)\cdot \sum_j \kj \lj \ps_j (v_j)^2.
\end{equation}
This implies that $\left(\phi(p^t) + \sum_j \kj \lj \ps_j (v_j)^2\right)$ decreases linearly, and finishes the proof of Theorem \ref{thm:tat-with-warehouse}, except that we need to show Conditions B3 and B4 hold throughout the tatonnement process.

\newpage

\subsection{Warehouse Stocks Are Bounded}\label{sect:bounding-warehouse}

So far we need $\kj$ to satisfy Conditions B2, B3 and B4.
Conditions B2 is satisfied so long as $\kj$ is sufficiently small.
However, we need to be cautious with Conditions B3 and B4 as it is not immediately evident
that $v_j$ remains bounded throughout the tatonnement process.

We begin with Phase 1.
The initial value of $\Xi_1$ decreases as $\kj$ decreases,
and Phase 1 ends when $\Xi_1$ is smaller than $\psi/2$, which is independent of $\kj$.
By \eqref{eq:Xi1-cont-concrete}, $\Xi_1$ drops linearly at a rate that does not depend on $\kj$.
Hence, the length of Phase 1 is finitely bounded when the $\kj$ are sufficiently small.
The change to each warehouse $j$ is upper bounded by
$$\mbox{(The length of Phase 1)}~\times ~\mbox{(Maximum excess demand for good $j$ in Phase 1)},$$
which is also finitely bounded.
This allows us to set $\kj$ sufficiently small to ensure that Conditions B3 and B4 hold throughout Phase 1.

Next, we consider Phase 2, which starts at some time $t_2$.
At $t_2$, which is the finishing time of Phase 1, Conditions B3 and B4 hold.
Let $B := \Xi_2(t_2)$.
Note that by \eqref{eq:Xi2-lower-concrete},
when Conditions B1--B4 hold,
there exist constants $C_1,C_2$ such that
\begin{equation}\label{eq:Xi2-lower-concrete-rewrite}
\Xi_2(t) \geq C_1 \phi(p^t) + C_2 \sum_j \kj \lj \ps_j (v_j)^2.
\end{equation}
We impose two additional conditions on $\kj$:
\begin{enumerate}
\item[B5.] $\kj$ are sufficiently small such that for all $j$, $\kj \leq \frac{C_2 \ps_j \lj}{101 B}$.
\item[B6.] $\kj$ are sufficiently small such that for all $j$, $\kj \leq \frac{C_2 \psi}{2(26/D_1+4) B}$.
\end{enumerate}

Suppose that at some time $t_3 > t_2$, Condition B3 or B4 ceases to hold.
By our analysis of Phase 2, $\Xi_2$ decreases between times $t_2$ and $t_3$, so $\Xi_2(t_3) \leq B$.

If Condition B3 ceases to hold at $t_3$, as the warehouse contents change smoothly, there exists a good $\ell$ with $|\kappa_\ell v_\ell| = 1/10$, and for other goods Condition B3 remains valid.
Thus we can still apply \eqref{eq:Xi2-lower-concrete-rewrite} with Condition B5 to yield
$$\Xi_2(t_3)\geq C_2 \kappa_\ell \lambda_\ell \ps_\ell (v_\ell)^2 = \frac{C_2 \lambda_\ell \ps_\ell}{\kappa_\ell} |\kappa_\ell v_\ell|^2 = \frac{C_2 \lambda_\ell \ps_\ell}{100 \kappa_\ell} > B,$$
which is a contradiction.

If Condition B4 ceases to hold at $t_3$, as the warehouse contents change smoothly, $\sum_j p_j \lj (\kv)^2 = \frac{1}{26/D_1+4}\psi$.
Thus we can still apply \eqref{eq:Xi2-lower-concrete-rewrite} with Condition B6 to yield
$$\Xi_2(t_3)\geq C_2 \sum_j \kj \lj \ps_j (v_j)^2 \geq \frac{C_2}{1.9 \max_j \kj} \sum_j p_j \lj (\kv)^2 \geq \frac{C_2}{1.9\kappa_j} \cdot \frac{1}{26/D_1+4}\psi > B,$$
which is a contradiction.

Thus, there does not exist $t_3>t_2$ at which Condition B3 or B4 ceases to hold,
i.e.~the two conditions hold throughout Phase 1 and Phase 2.

\newpage

\subsection{Missing Proofs}\label{sect:ongoing-missing-proof}

\begin{pfof}{Lemma \ref{lem:progress-of-phi-with-warehouse}}
We start with the proof of \eqref{eq:progress-of-phi-1}.
By Result (3) about Complementary CES markets (see the appendix on tatonnment):
\begin{align}
\nonumber
& \phi^- - \phi^+ ~\geq ~ [\tz_j + (z_j - \tz_j)] (\Dp_j) - \frac{1.5 x_j}{p_j} (\Dp_j)^2\\
\label{eq:bnd-for-next-deriv}
&\geq \tz_j (\Dp_j) - \frac{1.5 x_j}{p_j} (\Dp_j)^2 - |z_j - \tz_j|\cdot |\Dp_j|\\
\nonumber
&= \tz_j \frac{(\tz_j - \kv) \Dt_j}{\gjt} - \frac{1.5 x_j}{p_j}\left(\frac{(\tz_j-\kv) \Dt_j}{\gjt}\right)^2 - |z_j - \tz_j|\cdot |\Dp_j|\\
\nonumber
&\geq \tz_j \frac{(\tz_j - \kv) \Dt_j}{\gjt} - \frac{1}{2} \frac{(\tz_j-\kv)^2 \Dt_j}{\gjt} - |z_j - \tz_j|\cdot |\Dp_j|\comm{By Eqn.~\eqref{eq:progress-of-phi-inter} and $\Dt_j\leq 1$}\\
\nonumber
&= \frac{1}{2} \frac{(\tz_j)^2 \Dt_j}{\gjt} - \frac{1}{2} \frac{(\kv)^2 \Dt_j}{\gjt} - |z_j - \tz_j|\cdot |\Dp_j|.
\end{align}
\noindent
Next, we give the proof of \eqref{eq:progress-of-phi-2}.
From~\eqref{eq:bnd-for-next-deriv}:
\begin{align*}
& \phi^- - \phi^+ ~\geq ~ (\tz_j-\kv) (\Dp_j) - \frac{1.5 x_j}{p_j} (\Dp_j)^2 - |z_j - \tz_j|\cdot |\Dp_j| - |\kv|\cdot |\Dp_j|\\
&\geq \frac{\gjt \Dp_j}{\Dt_j} \cdot \Dp_j - \frac{1.5 x_j}{p_j} \frac{1}{\gjt}\cdot \gjt (\Dp_j)^2
- |z_j - \tz_j|\cdot |\Dp_j| - \frac{1}{2}\left(2\frac{(\kv)^2 \Dt_j}{\gjt} + \frac{1}{2} \frac{\gjt (\Dp_j)^2}{\Dt_j}\right)\\
& \hspace*{1.7in}\comm{For the last term use the AM-GM ineq.}\\
&\geq \frac{\gjt (\Dp_j)^2}{\Dt_j} - \frac{1}{15} \frac{\gjt (\Dp_j)^2}{\Dt_j} - \frac{(\kv)^2 \Dt_j}{\gjt} - \frac{1}{4} \frac{\gjt (\Dp_j)^2}{\Dt_j} - |z_j - \tz_j|\cdot |\Dp_j|\\
& \hspace*{1.7in}\comm{For the second term use Eqn.~\eqref{eq:progress-of-phi-inter} and $\Dt_j\leq 1$}\\
&= \frac{41}{60} \frac{\gjt (\Dp_j)^2}{\Dt_j} - \frac{(\kv)^2 \Dt_j}{\gjt} - |z_j - \tz_j|\cdot |\Dp_j|.
\end{align*}
\end{pfof}

\begin{pfof}{Equation \eqref{eq:Xi1-cont}}
Note that $\frac{d v_j}{dt} = -z_j^t$.
\begin{align*}
& \frac{d \Xi_1}{dt} \\
& =  -c_1 \sum_j \frac{(z_j^t)^2}{\gjsj} - c_2 \sum_j \sum_i
\xi^{\beta_i}_j \UH{k_i}{j}{\beta_i}{\sigj}{p_j^{\tau_j+}}\frac{(\Dp_{k_i})^2}{\Dt_{k_i}}
+ \sum_j \frac{(\kj v_j^t)^2}{\gjsj} - 2 \sum_j \frac{(\kj)^2 v_j^t z_j^t (t - \tau_j)}{\gjsj}\\
& \leq  -c_1 \sum_j \frac{(z_j^t)^2}{\gjsj} - c_2 \sum_j \sum_i
\xi^{\beta_i}_j \UH{k_i}{j}{\beta_i}{\sigj}{p_j^{\tau_j+}}\frac{(\Dp_{k_i})^2}{\Dt_{k_i}}
+ \sum_j \frac{(\kj v_j^t)^2}{\gjsj} + 2 \sum_j \frac{\kj}{\gjsj} \left|\kj v_j^t\right|\cdot |z_j^t|\\
& \leq  -c_1 \sum_j \frac{(z_j^t)^2}{\gjsj} - c_2 \sum_j \sum_i
\xi^{\beta_i}_j \UH{k_i}{j}{\beta_i}{\sigj}{p_j^{\tau_j+}}\frac{(\Dp_{k_i})^2}{\Dt_{k_i}}
+ \sum_j \frac{(\kj v_j^t)^2}{\gjsj} + \sum_j \frac{\kj}{\gjsj} \left[(\kj v_j^t)^2 + (z_j^t)^2\right]\\
& \hspace*{2.7in}\comm{For the last term use the AM-GM ineq.}\\
& =  -\sum_j (c_1-\kj) \frac{(z_j^t)^2}{\gjsj} + \sum_j (1+\kj) \frac{(\kj v_j^t)^2}{\gjsj} - c_2 \sum_j \sum_i
\xi^{\beta_i}_j \UH{k_i}{j}{\beta_i}{\sigj}{p_j^{\tau_j+}}\frac{(\Dp_{k_i})^2}{\Dt_{k_i}}.
\end{align*}
\end{pfof}

\begin{pfof}{Lemma \ref{lem:Xi1-at-update}}
\begin{align*}
\Xi_1^- - \Xi_1^+ &= \phi^- - \phi^+ - c_1 \int_{\tau_j}^t \frac{(z_j(t'))^2}{\gjt}\,dt' + \sum_i \xi^{\beta_i}_j
\UH{k_i}{j}{\beta_i}{t}{p_j^{\tau_j+}} \frac{(\Dp_{k_i})^2}{\Dt_{k_i}} [2-c_2(t-\beta_i)]\\
&\qquad\qquad - 2 \sum_{k\neq j} \xi^t_k \cdot \UH{j}{k}{t}{\sigk}{p_k^{\tau_k+}}\frac{(\Dp_j)^2}{\Dt_j}
+ \frac{(\kv)^2 \Dt_j}{\gjt}\\
&\geq \frac{1}{2}\left(\frac{41}{60} \frac{\gjt (\Dp_j)^2}{\Dt_j} - \frac{(\kv)^2 \Dt_j}{\gjt} - |z_j - \tz_j|\cdot |\Dp_j|\right)
~~~~~~~~~~~~\comm{By Eqn.~\eqref{eq:progress-of-phi-2}}\\
&\qquad\qquad + \frac{1}{2}\left(\frac{1}{2} \frac{(\tz_j)^2 \Dt_j}{\gjt} - \frac{1}{2} \frac{(\kv)^2 \Dt_j}{\gjt} - |z_j - \tz_j|\cdot |\Dp_j|\right)
\comm{By Eqn.~\eqref{eq:progress-of-phi-1}}\\
&\qquad\qquad - c_1 \int_{\tau_j}^t \frac{(z_j(t'))^2}{\gjt}\,dt' + (2-c_2) \sum_{i=1}^m \xi^{\beta_i}_j \UH{k_i}{j}{\beta_i}{t}{p_j^{\tau_j+}} \frac{(\Dp_{k_i})^2}{\Dt_{k_i}}\\
&\qquad\qquad - 2 \sum_{k\neq j} \xi^t_k \cdot \UH{j}{k}{t}{\sigk}{p_k^{\tau_k+}}\frac{(\Dp_j)^2}{\Dt_j} + \frac{(\kv)^2 \Dt_j}{\gjt}\\
&\geq \frac{41}{120} \frac{\gjt (\Dp_j)^2}{\Dt_j} + \frac{1}{4} \frac{(\tz_j)^2 \Dt_j}{\gjt} - \underbrace{|z_j - \tz_j|\cdot |\Dp_j|}_{F_1}
- \underbrace{c_1 \int_{\tau_j}^t \frac{(z_j(t'))^2}{\gjt}\,dt'}_{F_2}\\
&\qquad\qquad + (2-c_2) \sum_{i=1}^m \xi^{\beta_i}_j \UH{k_i}{j}{\beta_i}{t}{p_j^{\tau_j+}} \frac{(\Dp_{k_i})^2}{\Dt_{k_i}}
- \underbrace{2 \sum_{k\neq j} \xi^t_k \cdot \UH{j}{k}{t}{\sigk}{p_k^{\tau_k+}}\frac{(\Dp_j)^2}{\Dt_j}}_{F_3}.
\end{align*}
Note that $F_1$, $F_2$ and $F_3$ are similar to
the terms $E_1$, $E_2$ and $E_3$ in the proof of Lemma \ref{lem:Phi-at-update}.
We can bound $F_1,F_2,F_3$ similarly to the way we bounded $E_1,E_2,E_3$.

Recall from the proof of Lemma \ref{lem:Phi-at-update} that $V_2 := \sum_{i=1}^m \xi^{\beta_i}_j \UH{k_i}{j}{\beta_i}{t}{p_j^{\tau_j+}} \frac{(\Dp_{k_i})^2}{\Dt_{k_i}}$. We derive the following bounds:
\begin{align*}
F_1 & \leq 2\eB \gjt (\Dp_j)^2 + V_2;\\
F_2 & \leq c_1 (1+4\eB) \frac{(\tz_j)^2 \Dt_j}{\gjt} + c_1(2+8\eB) V_2;\\
F_3 & \leq 2\eF \gjt \frac{(\Dp_j)^2}{\Dt_j}.
\end{align*}

Thus
\begin{align*}
\Xi_1^- - \Xi_1^+ &\geq \left(\frac{41}{120} - 2 \eB - 2\eF\right) \frac{\gjt(\Dp_j)^2}{\Dt_j} + \left(\frac{1}{4} - c_1(1+4\eB)\right)\frac{(\tz_j)^2 \Dt_j}{\gjt}\\
&\qquad + \left(1-c_2 - c_1(2+8\eB)\right) V_2.\end{align*}
Note that by
Eqns.~\eqref{eq:eF-bound}
and~\eqref{eq:eB-bound},
$2\eF + 2\eB = 0.304 < \frac{41}{120}$, $1+4\eB < 1.4$ and $2+8\eB < 2.7$.
The result now follows.
\end{pfof}

\begin{pfof}{Lemma \ref{lem:Phi-at-update-with-warehouse}}
This proof is similar to the one of Lemma \ref{lem:Xi1-at-update}; we only point out the key steps.
\begin{align*}
\Phi^- - \Phi^+ &\geq \phi^- - \phi^+ - c_1 \int_{\tau_j}^t \lj p_j (z_j(t'))^2\,dt'
+ (6-c_2)\sum_i \xi^{\beta_i}_j \UH{k_i}{j}{\beta_i}{t}{p_j^{\tau_j+}} \frac{(\Dp_{k_i})^2}{\Dt_{k_i}}\\
&\qquad\qquad - 6 \sum_{k\neq j} \xi^t_k \cdot \UH{j}{k}{t}{\sigk}{p_k^{\tau_k+}}\frac{(\Dp_j)^2}{\Dt_j}\\
&\geq \frac{9}{10}\left(\frac{41}{60} \frac{(\Dp_j)^2}{\lj p_j \Dt_j} - \lj p_j (\kv)^2 \Dt_j - |z_j - \tz_j|\cdot |\Dp_j|\right)
\comm{By Eqn.~\eqref{eq:progress-of-phi-2}}\\
&\qquad\qquad + \frac{1}{10}\left(\frac{1}{2} \lj p_j (\tz_j)^2 \Dt_j - \frac{1}{2} \lj p_j (\kv)^2 \Dt_j - |z_j - \tz_j|\cdot |\Dp_j|\right)
\comm{By Eqn.~\eqref{eq:progress-of-phi-1}}\\
&\qquad\qquad - c_1 \int_{\tau_j}^t \lj p_j (z_j(t'))^2\,dt' + (6-c_2) \sum_{i=1}^m \xi^{\beta_i}_j \UH{k_i}{j}{\beta_i}{t}{p_j^{\tau_j+}} \frac{(\Dp_{k_i})^2}{\Dt_{k_i}}\\
&\qquad\qquad - 6 \sum_{k\neq j} \xi^t_k \cdot \UH{j}{k}{t}{\sigk}{p_k^{\tau_k+}}\frac{(\Dp_j)^2}{\Dt_j}\\
&\geq \frac{123}{200} \frac{(\Dp_j)^2}{\lj p_j \Dt_j} + \frac{1}{20} \lj p_j (\tz_j)^2 \Dt_j - \frac{19}{20} \lj p_j (\kv)^2 \Dt_j - \underbrace{|z_j - \tz_j|\cdot |\Dp_j|}_{F_1} \\
&\qquad\qquad - \underbrace{c_1 \int_{\tau_j}^t \lj p_j (z_j(t'))^2\,dt'}_{F_2} + (6-c_2) \sum_{i=1}^m \xi^{\beta_i}_j \UH{k_i}{j}{\beta_i}{t}{p_j^{\tau_j+}} \frac{(\Dp_{k_i})^2}{\Dt_{k_i}}\\
&\qquad\qquad - \underbrace{6 \sum_{k\neq j} \xi^t_k \cdot \UH{j}{k}{t}{\sigk}{p_k^{\tau_k+}}\frac{(\Dp_j)^2}{\Dt_j}}_{F_3 '}.
\end{align*}

Then we apply the bounds on $F_1,F_2,F_3$ in the proof of Lemma \ref{lem:Xi1-at-update}.\footnote{There
is one minor difference: $\gamma_j^{\sigj}$ is replaced by $1/(\lj p_j)$.
Also, $F_3 '$ is three times the value of $F_3$, so the bound on $F_3 '$ is amplified accordingly.}
to show that
\begin{align*}
\Phi^- - \Phi^+ &\geq \left(\frac{123}{200} - 2 \eB - 6\eF\right) \frac{(\Dp_j)^2}{\lj p_j \Dt_j} + \left(\frac{1}{20} - c_1(1+4\eB)\right) \lj p_j (\tz_j)^2 \Dt_j\\
&\qquad - \frac{19}{20} \lj p_j (\kv)^2 \Dt_j + \left(5 - c_2 - c_1(2+8\eB)\right) V_2.\end{align*}
Note that $\frac{123}{200} - 2\eB - 6\eF = 0.039$, $1+4\eB < 1.4$ and $2+8\eB < 2.7$; the lemma now follows.
\end{pfof}

\begin{pfof}{Equation~\eqref{eq:phase2-W-cont}}
Note that $\frac{d v_j}{dt} = -z_j^t$.
\begin{align*}
\frac{d\wh}{dt} &= \sum_j \ps_j \left[\frac{2\kj}{\lj}(f_j + \lj v_j^t)(-\lj z_j^t) - c_3 \lj (\kv^t)^2 + 2 c_3 \lj (\kj)^2 v_j^t z_j^t (t-\tau_j) + 2\kj \lj v_j^t z_j^t\right]\\
&\leq \sum_j \ps_j \left[2\kj |f_j| |z_j^t| - c_3 \lj (\kv^t)^2 + 2 c_3 \lj \kj |\kv^t| |z_j^t|\right]\\
&\leq \sum_j \ps_j \left[\kj\left(\frac{(f_j)^2}{R_1} + R_1 (z_j^t)^2\right) - c_3 \lj (\kv^t)^2 + c_3 \lj \kj \left[(\kv^t)^2 + (z_j^t)^2\right]\right]\\
&= -c_3 \sum_j (1-\kappa_j) \lj \ps_j (\kv^t)^2 + \sum_j (R_1 + c_3 \lj) \kj \ps_j (z_j^t)^2  + \frac{1}{R_1}\sum_j \kj \ps_j (f_j)^2.
\end{align*}
\end{pfof}

\begin{pfof}{Lemma \ref{lem:W-at-update-with-warehouse}}
At the price update, $f_j^+ = f_j^- + \ln\left(1 + \lj (\tz_j - \kv) \Dt_j\right)$.
Note that in Phase 2, $\left|\lj (\tz_j - \kv) \Dt_j \right| \leq 1/60$
and hence $\ln\left(1 + \lj (\tz_j - \kv) \Dt_j\right) = (1+\chi) \lj (\tz_j - \kv) \Dt_j$ for some $\chi$ with $|\chi | \leq \frac{1}{100}$.\footnote{When
$|y| \leq \frac{1}{60}$, $\ln(1+y) \in \left[1-\frac{1}{100},1+\frac{1}{100}\right]\cdot y$.}
Then
\begin{align*}
\wh^- - \wh^+ &= \ps_j\left[\frac{\kj}{\lj} \left[(f_j + \lj v_j)^2 - (f_j + (1+\chi) \lj (\tz_j - \kv) \Dt_j + \lj v_j)^2\right] \right.\\
&\qquad\qquad \left. - c_3 \lj (\kv)^2 \Dt_j + 2 \kj \lj \int_{\tau_j}^t v_j(t') z_j(t')\,dt' \right]
\end{align*}

Let $\bz_j$ be the average excess demand for good $j$ between times $\tau_j$ and $t$, i.e.~$\bz_j := \frac{1}{t_2-t_1} \int_{t_1}^{t_2} z_j^{t'}\,dt'$.
Note that $v_j(\tau_j) = v_j(t) + \bz_j \Dt_j$ and $\frac{d v_j}{d t} = -z_j$. We use integration by substitution to evaluate the integral in the above formula:
$$\int_{\tau_j}^t v_j(t') z_j(t')\,dt' = -\int_{v_j(\tau_j)}^{v_j(t)} v_j\,dv_j = \frac{1}{2}\left(v_j(\tau_j)^2 - v_j(t)^2\right) = v_j \bz_j \Dt_j + \frac{1}{2} (\bz_j)^2 (\Dt_j)^2.$$
By direct expansion and regrouping terms, we have
\begin{align*}
&\wh^- - \wh^+ \\
&= \ps_j \Dt_j \left\{ \left[2(1+\chi) - (1+\chi)^2 \kj \Dt_j - c_3\right] \lj (\kv)^2 + \left[(\bz_j)^2 - (\tz_j)^2\right] \kj \lj \Dt_j\right.\\
&\qquad\qquad - (2\chi + \chi^2) \kj \lj (\tz_j)^2 \Dt_j + 2\lj (\bz_j - \tz_j) \kv \\
&\qquad\qquad + \left.\left[(1+\chi)^2 \kj \Dt_j - \chi\right] \cdot 2\lj \tz_j \kv - 2(1+\chi) \kj f_j (\tz_j - \kv)\right\} \\
&\geq \ps_j \Dt_j \left\{\left[2(1+\chi) - (1+\chi)^2 \kj \Dt_j - c_3\right] \lj (\kv)^2 - \underbrace{\kj \lj |(\bz_j)^2 - (\tz_j)^2|}_{G_1}\right.\\
&\qquad\qquad - \underbrace{|2\chi + \chi^2|\cdot \kj \lj (\tz_j)^2}_{G_2} - \underbrace{2 |\bz_j - \tz_j|\cdot \left|\lj \kv\right|}_{G_3}\\
&\qquad\qquad \left.- \underbrace{2 |(1+\chi)^2 \kj \Dt_j - \chi| \lj |\tz_j| \cdot |\kv|}_{G_4} - \underbrace{2(1+\chi) \kj |f_j| \cdot |\tz_j - \kv|}_{G_5}\right\}
\end{align*}
Next, we bound the terms $G_1,G_2,G_3,G_4,G_5$.
Recall the notations we use in the proof of Lemma \ref{lem:Phi-at-update}
$V_1 := \sum_{k\neq j} \frac{1}{\min_{i:k_i=k} \xi^{\beta_i}_j} \UH{k}{j}{\tau_j}{t}{p_j^t}$ and
$V_2 := \sum_{i=1}^m \xi^{\beta_i}_j\cdot \UH{k_i}{j}{\beta_i}{t}{p_j^t} \frac{(\Dp_{k_i})^2}{\Dt_{k_i}}.$

\begin{align*}
G_1 &\leq \kj \lj \left[ (\bz_j - \tz_j)^2 + \frac{2}{\lj p_j} |\lj p_j \tz_j| \cdot |\bz_j - \tz_j| \right]\\
&\leq \kj \lj \left[ 8 V_1 V_2 + \frac{2}{\lj p_j}\left(2(\lj p_j \tz_j)^2 V_1 + V_2\right)\right] \comm{By Eqns.~\eqref{eq:error-of-gradient-3} and \eqref{eq:error-of-gradient-2}}\\
&\leq \kj \lj \left[ 8 \frac{\eB}{\lj p_j} V_2 + 4 \eB (\tz_j)^2 + \frac{2}{\lj p_j} V_2\right]
\comm{as by Cond.\ (A2), $V_1 \le \eB\gamma_j^t = \eB/(\lambda_j p_j)$}\\
&= 4 \eB \kj \lj (\tz_j)^2 + \frac{(2+8\eB)\kj}{p_j} V_2.
\end{align*}

To bound $G_2$, note that $|\chi| \leq 1/100$ and $\kj \leq 1/600$ imply that
$\kj |2\chi + \chi^2| \leq 0.0000335$,
and hence $G_2 \leq 0.0000335 \lj (\tz_j)^2$.

\begin{align*}
G_3 &= \frac{2}{p_j}|\bz_j - \tz_j| \cdot |\lj p_j \kv|\\
&\leq \frac{2}{p_j}\left[2 (\lj p_j \kj v_j)^2 V_1 + V_2\right]\comm{By Eqn.~\eqref{eq:error-of-gradient-2}}\\
&\leq \frac{2}{p_j}\left[2 (\lj p_j \kj v_j)^2 \frac{\eB}{\lj p_j} + V_2\right]\\
&= 4\eB \lj (\kv)^2 + \frac{2}{p_j} V_2.
\end{align*}

To bound $G_4$, note that $|\chi|\leq 1/100$, $\kj \leq 1/600$ and $\Dt_j\leq 1$ imply that
$|(1+\chi)^2 \kj \Dt_j - \chi| \leq 0.0117$.
Then by AM-GM inequality, $G_4 \leq 0.0117 \lj (\tz_j)^2 + 0.0117 \lj (\kv)^2$.

\begin{align*}
G_5 &= 2(1+\chi) \frac{\kj}{\lj p_j} |f_j|\cdot |\lj p_j (\tz_j - \kv)| \\
&= 2(1+\chi) \frac{\kj}{\lj p_j} |f_j| \cdot \left| \frac{\Dp_j}{\Dt_j} \right|\\
&\leq \frac{101}{100} \frac{\kj}{\lj p_j} \left(\frac{\kj p_j (f_j)^2}{R_2} + \frac{R_2 (\Dp_j)^2}{\kj p_j (\Dt_j)^2}\right).\comm{by the AM-GM ineq.}
\end{align*}

Combining all the above bounds yields
\begin{align*}
\wh^- - \wh^+ &\geq \left[2(1+\chi) - (1+\chi)^2 \kj \Dt_j - c_3 - 4\eB - 0.0117\right]\frac{\ps_j}{p_j} \lj p_j (\kv)^2 \Dt_j\\
&\qquad -(0.0118 + 4 \eB \kj)\frac{\ps_j}{p_j} \lj p_j (\tz_j)^2 \Dt_j - \frac{(2+2\kj + 8\eB \kj) \ps_j}{p_j} V_2\\
&\qquad - \frac{101}{100}\cdot \frac{\kj}{\lj}\cdot \frac{\ps_j}{p_j} \left(\frac{\kj p_j (f_j)^2 \Dt_j}{R_2} + \frac{R_2 (\Dp_j)^2}{\kj p_j \Dt_j}\right).
\end{align*}

Note the following:
\begin{itemize}
\item $|\chi|\leq 1/100$, $\kj \leq 1/600$ and $\Dt_j\leq 1$ imply that
$2(1+\chi) - (1+\chi)^2 \kj \Dt_j \geq 1.9783$. Also, recall that $\eB = 0.084$.
Thus $\left[2(1+\chi) - (1+\chi)^2 \kj \Dt_j - c_3 - 4\eB - 0.0117\right]\frac{\ps_j}{p_j} \geq (1.6306-c_3)/1.9 \geq 0.858 - c_3/1.9$.
\item $\eB = 0.084$ and $\kj\leq 1/600$ imply that $(0.0118 + 4\eB \kj)\frac{\ps_j}{p_j} \leq 0.01236\times 1.9 \leq 0.0235$.
\item $\eB = 0.084$ and $\kj\leq 1/600$ imply that $\frac{(2+2\kj + 8\eB \kj) \ps_j}{p_j} \leq 2.00446\times 1.9 \leq 3.809$.
\item $\frac{101}{100}\cdot \frac{\ps_j}{p_j} \leq 1.92$.
\end{itemize}
The lemma follows.
\end{pfof}

In the proofs of Equations \eqref{eq:Xi2-lower-bound} and \eqref{eq:Xi2-upper-bound} below,
we need the following bound on $(\bz_j)^2$:
\begin{align*}
(\bz_j)^2 - (z_j)^2 &= (\bz_j - z_j)^2 - 2 z_j (z_j - \bz_j)\\
&\leq 8 V_1 V_2 + \frac{1}{5 \lj p_j} \left| 10 \lj p_j z_j\right|\cdot |z_j - \bz_j|
~~~~~~~~~\comm{by Eqn.~\ref{eq:error-of-gradient-3}}\\
&\leq \frac{8\eB}{\lj p_j} V_2 + \frac{1}{5 \lj p_j}\left(200 (\lj p_j)^2 (z_j)^2 V_1 + V_2\right)
\comm{as $V_1 \le \eB/(\lambda_j p_j)$}\\
&\leq \frac{0.672}{\lj p_j} V_2 + 40 \lj p_j (z_j)^2 \frac{\eB}{\lj p_j} + \frac{0.2}{\lj p_j} V_2\\
&= 3.36 (z_j)^2 + \frac{0.872}{\lj p_j} \sum_i \xi^{\beta_i}_j \UH{k_i}{j}{\beta_i}{\sigj}{p_j^{\tau_j+}} \frac{(\Dp_{k_i})^2}{\Dt_{k_i}}
\end{align*}
and hence
\begin{equation}
\label{eq:pj-zjbar-bound}
\lj p_j (\bz_j)^2 \leq 4.36 \lj p_j (z_j)^2 + 0.872 \sum_i \xi^{\beta_i}_j \UH{k_i}{j}{\beta_i}{\sigj}{p_j^{\tau_j+}} \frac{(\Dp_{k_i})^2}{\Dt_{k_i}}.
\end{equation}

\begin{pfof}{Equation \eqref{eq:Xi2-lower-bound}}
By Lemma \ref{lem:relating-Phi-phi}, if $2-c_2 \geq 2.7 c_1$,
then
$$\phi(p^t) - c_1 \sum_j \int_{\tau_j}^{t} \lj p_j (z_j(t'))^2\,dt' + \sum_j \sum_i \xi^{\beta_i}_{j} \UH{k_i}{j}{\beta_i}{\sigj}{p_j^{\tau_j+}} \frac{(\Dp_{k_i})^2}{\Dt_{k_i}} \left[2-c_2 (t - \beta_i)\right] \geq (1- 2.7 c_1) \phi(p^t).$$
Thus, $\Phi$, as defined in \ref{eq:Xi2-Phi}, satisfy
\begin{align*}
\Phi &= \phi(p^t) - c_1 \sum_j \int_{\tau_j}^{t} \lj p_j (z_j(t'))^2\,dt' + \sum_j \sum_i \xi^{\beta_i}_{j} \UH{k_i}{j}{\beta_i}{\sigj}{p_j^{\tau_j+}} \frac{(\Dp_{k_i})^2}{\Dt_{k_i}} \left[6-c_2 (t - \beta_i)\right]\\
&\geq (1- 2.7 c_1) \phi(p^t) + 4 \sum_j \sum_i \xi^{\beta_i}_j \UH{k_i}{j}{\beta_i}{\sigj}{p_j^{\tau_j+}} \frac{(\Dp_{k_i})^2}{\Dt_{k_i}}.
\end{align*}

\begin{align*}
& \wh \\
&= \sum_j \frac{\kj}{\lj} \ps_j \left(f_j + \lj v_j\right)^2 - c_3 \sum_j \lj \ps_j (\kv)^2 (t-\tau_j)  + 2 \sum_j \kj \lj \ps_j \int_{\tau_j}^t v_j(t') z_j(t')\,dt'\\
&\geq \sum_j \frac{\kj}{\lj} \ps_j \left( \frac{(\lj v_j)^2}{2} - (f_j)^2 \right)
-c_3 \sum_j \lj \ps_j (\kv)^2
+ 2 \sum_j \kj \lj \ps_j \left(v_j \bz_j (t-\tau_j) + \frac{1}{2} (\bz_j)^2 (t-\tau_j)^2\right)\\
&\geq \sum_j \left(\frac{1}{2} - c_3 \kj\right) \lj \kj \ps_j (v_j)^2
- \sum_j \frac{\kj}{\lj} \ps_j (f_j)^2\\
&\qquad\qquad + 2 \sum_j \kj \lj \ps_j \left(-\frac{1}{6} (v_j)^2 - \frac{3}{2} (\bz_j)^2 (t-\tau_j)^2 + \frac{1}{2} (\bz_j)^2 (t-\tau_j)^2\right)
\comm{by the AM-GM ineq.}\\
&\geq \sum_j \left(\frac{1}{6} - c_3 \kj\right) \lj \kj \ps_j (v_j)^2
- \sum_j \frac{\kj}{\lj} \ps_j (f_j)^2
- 2\sum_j \kj \lj \ps_j (\bz_j)^2 \\
&\geq \sum_j \left(\frac{1}{6} - c_3 \kj\right) \lj \kj \ps_j (v_j)^2
- \sum_j \frac{\kj}{\lj} \ps_j (f_j)^2\\
&\qquad\qquad - 3.8 \sum_j \kj \left(4.36 \lj p_j (z_j)^2 + 0.872 \sum_i \xi^{\beta_i}_j \UH{k_i}{j}{\beta_i}{\sigj}{p_j^{\tau_j+}} \frac{(\Dp_{k_i})^2}{\Dt_{k_i}}\right)
\comm{by eqn.~\eqref{eq:pj-zjbar-bound}}\\
&\geq \sum_j \left(\frac{1}{6} - c_3 \kj\right) \lj \kj \ps_j (v_j)^2
- \sum_j \frac{\kj}{\lj} \ps_j (f_j)^2\\
&\qquad\qquad - 16.6 \sum_j \kj \lj p_j (z_j)^2
- 3.314 \sum_i \xi^{\beta_i}_j \UH{k_i}{j}{\beta_i}{\sigj}{p_j^{\tau_j+}} \frac{(\Dp_{k_i})^2}{\Dt_{k_i}}.
\end{align*}
Recall that $\Xi_2 = \Phi + 1.2 \wh + 0.1212 \sum_j \frac{\kj \ps_j (f_j)^2}{R_2} (t-\tau_j) \geq \Phi + 1.2 \wh$. With the two inequalities above, the result follows.
\end{pfof}

\begin{pfof}{Equation \eqref{eq:Xi2-upper-bound}}
It follows immediately from~\eqref{eq:Xi2-Phi} that
$$\Phi \leq \phi(p^t) + 6 \sum_j \sum_i \xi^{\beta_i}_j \UH{k_i}{j}{\beta_i}{\sigj}{p_j^{\tau_j+}} \frac{(\Dp_{k_i})^2}{\Dt_{k_i}}.$$

\begin{align*}
\wh &= \sum_j \frac{\kj}{\lj} \ps_j \left(f_j + \lj v_j\right)^2 - c_3 \sum_j \lj \ps_j (\kv)^2 (t-\tau_j)  + 2 \sum_j \kj \lj \ps_j \int_{\tau_j}^t v_j(t') z_j(t')\,dt'\\
&\leq 2 \sum_j \frac{\kj}{\lj} \ps_j (f_j)^2 + 2 \sum_j \kj \lj \ps_j (v_j)^2
+ 2 \sum_j \kj \lj \ps_j \left(v_j \bz_j (t-\tau_j) + \frac{1}{2} (\bz_j)^2 (t-\tau_j)^2\right)\\
&\leq 2 \sum_j \frac{\kj}{\lj} \ps_j (f_j)^2 + 2 \sum_j \kj \lj \ps_j (v_j)^2
+ 2 \sum_j \kj \lj \ps_j \left(\frac{1}{2} (v_j)^2 + \frac{1}{2} (\bz_j)^2 (t-\tau_j)^2 + \frac{1}{2} (\bz_j)^2 (t-\tau_j)^2\right)\\
&\leq 2 \sum_j \frac{\kj}{\lj} \ps_j (f_j)^2 + 3 \sum_j \kj \lj \ps_j (v_j)^2
+ 2 \sum_j \kj \lj \ps_j (\bz_j)^2 \\
&\leq 2 \sum_j \frac{\kj}{\lj} \ps_j (f_j)^2 + 3 \sum_j \kj \lj \ps_j (v_j)^2\\
&\qquad\qquad + 3.8 \sum_j \kj \left(4.36 \lj p_j (z_j)^2 + 0.872 \sum_i \xi^{\beta_i}_j \UH{k_i}{j}{\beta_i}{\sigj}{p_j^{\tau_j+}} \frac{(\Dp_{k_i})^2}{\Dt_{k_i}}\right)
\comm{by eqn.~\eqref{eq:pj-zjbar-bound}}\\
&\leq 2 \sum_j \frac{\kj}{\lj} \ps_j (f_j)^2 + 3 \sum_j \kj \lj \ps_j (v_j)^2\\
&\qquad\qquad + 16.6 \sum_j \kj \lj p_j (z_j)^2
+ 3.314 \sum_i \xi^{\beta_i}_j \UH{k_i}{j}{\beta_i}{\sigj}{p_j^{\tau_j+}} \frac{(\Dp_{k_i})^2}{\Dt_{k_i}}.
\end{align*}
Recall that $\Xi_2 = \Phi + 1.2 \wh + 0.1212 \sum_j \frac{\kj \ps_j (f_j)^2}{R_2} (t-\tau_j) \leq \Phi + 1.2 \wh + 0.1212 \sum_j \frac{\kj \ps_j (f_j)^2}{R_2}$. With the two inequalities above, the result follows.
\end{pfof}

To prove Lemma \ref{lem:z-to-f}, we need the following lemma.

\begin{lemma}\label{lem:phi-convexity-bound}
For all $p'\in R(1.9)$, $\phi(p') \geq \frac{1-\bar{\theta}}{13.28}\sum_j p_j^*(f_j)^2.$
\end{lemma}
\begin{pf}
Let $x_{ij}(p')$ be the demand for good $j$ of buyer $i$ at price $p'$. Note that
$$\frac{\partial^2 \phi}{\partial (p_j)^2} (p') = \sum_i \left(\frac{\theta_i (x_{ij}(p'))^2}{e_i} + \frac{(1-\theta_i)x_{ij}(p')}{p'_j}\right)\,\,\,\,\mbox{and}\,\,\,\,\,\,\frac{\partial^2 \phi}{\partial p_j \partial p_k} (p') = \sum_i \frac{\theta_i x_{ij}(p')x_{ik}(p')}{e_i}.$$

Let $A^i(p')$ denote the matrix with $A^i_{jk}(p') = x_{ij}(p')x_{ik}(p')$. Let $B^i(p')$ denote the diagonal matrix with $B^i_{jj}(p') = x_{ij}(p')/p'_j$.
Then the Hessian of $\phi$ at $p'$, which we denote it by $H(p')$, is $\sum_i \frac{\theta_i}{e_i} A^i(p') + \sum_i (1-\theta_i) B^i(p')$.

There are two key observations: first that $A^i$ is positive semi-definite and second that $\sum_i (1-\theta_i) B^i(p')$ majorizes $(1-\bar{\theta})\sum_i B^i(p')$,
where $\bar{\theta} = \max_i \theta_i$.
Hence $H(p')$ majorizes $(1-\bar{\theta})\sum_i B^i(p') := (1-\bar{\theta}) B(p')$, where $B_{jj}(p') = x_j(p')/p'_j$.
As $p'\in R(1.9)$, $x_j(p')\geq 1/1.9$ and $p'_j\leq 1.9 p_j^*$. Hence $B_{jj}(p') \geq \frac{1}{3.61 p_j^*}$.

\newcommand{\bphi}{\bar{\phi}}

Next, consider the function $\bphi(p) = \phi(p) - \sum_j \frac{1-\bar{\theta}}{7.22 p_j^*} (p_j - p_j^*)^2$.
Observe that for all $j$, $\frac{\partial \bphi}{\partial p_j}(p^*) = 0$ and the Hessian of $\bphi$ at every $p'\in R$ majorizes the zero matrix;
consequently, $\bphi$ is convex in $R(1.9)$, and $p^*$ is its minimum point.
Note that $\bphi(p^*) = 0$, so for all $p'\in R(1.9)$, $\phi(p') \geq \sum_j\frac{1-\bar{\theta}}{7.22 p_j^*} (p'_j - p_j^*)^2$.

Since $\left(\frac{p'_j-p_j^*}{p_j^*}\right)^2 \geq 0.544 \ln^2 \frac{p_j'}{p_j^*} = 0.544 (f_j)^2$,
$\phi(p') - \phi^*\geq \frac{1-\bar{\theta}}{13.28}\sum_j p_j^*(f_j)^2$.
\end{pf}

\begin{pfof}{Lemma \ref{lem:z-to-f}}
\cite[Lemma 6.3]{CCD2013} showed that for all $p'\in R(1.9)$, $\phi(p') \leq \max\left\{2,\frac{\bar{\theta}}{2\left(1+\bar{\theta}-2^{\bar{\theta}}\right)}\right\}\cdot \sum_j p'_j(z_j)^2$. Combining this with Lemma \ref{lem:phi-convexity-bound} yields the result.
\end{pfof}

\end{document}